\documentclass[12pt]{amsart}
\usepackage{a4wide,amsfonts,graphicx,
amssymb,stmaryrd,amscd,amsmath,latexsym,amsbsy}
\numberwithin{equation}{section}
\theoremstyle{plain}
\newtheorem{thm}{Theorem}[section]

\newtheorem{prop}[thm]{Proposition}

\newtheorem{defi}[thm]{Definition}
\newtheorem{lem}[thm]{Lemma}
\newtheorem{cor}[thm]{Corollary}
\newtheorem{eg}[thm]{Example}

\theoremstyle{remark}
\newtheorem{rema}[thm]{Remark}

\title[Extended trigonometric Cherednik algebras]
{Extended trigonometric Cherednik algebras and nonstationary
Schr{\"o}dinger equations with delta-potentials}
\author{J.T. Hartwig \& J.V. Stokman}
\address{J.T.H.: Department of Mathematics, Stanford University, USA\\
J.V.S: Korteweg-de Vries 
Institute for Mathematics, Universiteit van Amsterdam,
Science Park 904, 1098 XH Amsterdam, The Netherlands}
\email{jonas.hartwig@gmail.com \& j.v.stokman@uva.nl}

\newcommand\al{\alpha}   
 \newcommand\De{\Delta}

   \newcommand\ph{\varphi}
   \newcommand\om{\omega} \newcommand\Om{\Omega}

\newcommand\hf{\mathfrak{h}} 
 \newcommand\gf{\mathfrak{g}}

\newcommand\C{\mathbb{C}} \newcommand\Z{\mathbb{Z}}

\newcommand{\wh}[1]{\widehat{#1}}

\newcommand{\Ra}{\widehat{R}}
\newcommand{\Rap}{\widehat{R}^+}
\newcommand{\Wa}{\widehat{W}}
\newcommand{\Fa}{\widehat{F}}

\newcommand\reg{\mathrm{reg}}

\newcommand{\Hom}{\mathrm{Hom}}

\newcommand{\GL}{\mathrm{GL}}

\begin{document}
\begin{abstract}
We realize an extended version of the trigonometric Cherednik algebra
as affine Dunkl operators involving Heaviside functions. We use the quadratic
Casimir element of the extended trigonometric Cherednik algebra to define
an explicit nonstationary Schr{\"o}\-din\-ger equation 
with delta-potential. We use 
coordinate Bethe ansatz methods to construct solutions of the
nonstationary Schr{\"o}dinger equation in terms of
generalized Bethe wave functions.
It is shown that 
the generalized Bethe wave functions satisfy affine difference 
Knizhnik-Zamolodchikov equations as functions of in the momenta. 
The relation to the vector valued root system analogs of the quantum
Bose gas on the circle with delta-function interactions is indicated.
\end{abstract}

\maketitle

\section{Introduction}

The one dimensional 
quantum Bose gas with pairwise
delta-function interactions \cite{LL} 
is among the first nontrivial quantum integrable systems 
that were successfully analyzed using
the coordinate Bethe ansatz \cite{LL,Y,YY}. On the circle the 
coordinate Bethe ansatz leads  
to Bethe wave functions, given as explicit plane wave expansions with 
momenta subject to Bethe ansatz equations, which solve the associated 
stationary Schr{\"o}dinger equation.
In this paper we will construct Bethe 
wave type functions solving nonstationary extensions of 
Schr{\"o}dinger equations associated to vector valued, root system analogs
of the quantum Bose gas on the circle with pairwise
delta-function interactions. 
We will show that the role of the Bethe ansatz equations is taken over by 
difference analogs of Knizhnik-Zamolodchikov equations.

It is well known that the quantum inverse scattering method can be applied
to the one dimensional quantum Bose gas with pairwise
delta-function interactions
(see \cite{KP} and references therein). The key point is the fact
that the pertinent quantum Bose gas arises as particle sector of the integrable
quantum field theory in $1+1$ dimensions governed by the 
nonlinear Schr{\"o}dinger equation. The resulting possibility to 
create and annihilate quantum particles in the quantum Bose
gas has in particular 
led to the explicit evaluation of the norms of the Bethe wave functions
(see \cite{K,KP}) in terms of the Hessian of the 
Yang-Yang action \cite{YY}. 

Quantum Calogero-Moser systems \cite{OPe,BFV}, which can
be defined for any root system, form an important class of one dimensional
integrable quantum many body systems with pairwise interactions of
rational, trigonometric, hyperbolic or elliptic type. The quantum trigonometric
Calogero-Moser systems are naturally related to harmonic analysis on 
symmetric spaces \cite{HS}. Substantial progress has been made
over the past decades in solving quantum Calogero-Moser systems, 
with main tool the explicit
realization of (degenerate) Hecke algebra symmetries in terms of Dunkl type 
differential-reflection operators (see, e.g., 
\cite{D,H,Op,Cherednik1997, Cherednik1995, BFV}).  
These systems do not arise though as
particle sectors in particular integrable quantum field theories, hence
the application of quantum inverse scattering techniques to such systems 
is limited. In particular there is no analog of particle
creation and annihilation, which for instance explains
the completely different techniques in deriving norm formulas; it is based on
shift operators or intertwiners in the Hecke algebra context
(see, e.g., \cite{Oshift, Cherednik1997, OpLecture}) and on quantum particle
creation/annihilation in the quantum inverse scattering context \cite{K}.

The one dimensional quantum Bose gas with pairwise
delta-function interactions,
which is completely accessible to quantum inverse scattering techniques, can 
also be naturally viewed as member of the family of integrable quantum 
Calogero-Moser type systems. Firstly, the one dimensional quantum Bose gas
with pairwise
delta-function interactions has natural root system generalizations,
going back to Gaudin, Gutkin and Sutherland \cite{Ga,GutSut1979,Gutkin1982}.
For classical root systems it corresponds to imposing
integrable reflecting boundary conditions on the quantum particles on the line.
Secondly, Dunkl type operators and integral-reflection operators have been
associated to the one dimensional quantum Bose gas with delta-function 
interactions and their root system generalizations, which opens the way to
apply the Hecke algebra techniques 
(see, e.g., \cite{Gutkin1982,Po,HO,EOS0,EOS}).

The quantum Bose gas with pairwise 
delta-function interactions thus is accessible
for an unusually large variety of techniques from integrable systems and 
representation theory. This feature places the one dimensional
quantum Bose gas with delta-function interactions center stage of
various new developments in mathematical 
physics and representation theory (see, e.g., \cite{HO} and \cite{GS}
for two striking examples) and leads to the intriguing question how 
the quantum inverse scattering method and the Hecke algebra method can 
be united. For instance, conjectures \cite{Em}
have been put forward extending
Korepin's \cite{K} norm formulas to Bethe wave functions of the
root system analogs of the quantum Bose gas with pairwise 
delta-function interactions on the circle. Since Korepin's arguments 
are no longer applicable in the context of arbitrary root systems
by the loss of quantum particle creation and 
annihilation techniques, it seems
to require now Hecke algebra techniques instead to properly understand such
norm formulas. 

This formed an important motivation for the research leading up to the
present work. Our starting point 
is the observation from \cite{EOS} that the root
system analogs of the quantum Bose gas on the circle with pairwise
delta-function
interactions arise from a representation of the trigonometric Cherednik
algebra {\it at critical level} in terms of Dunkl type operators.
This suggests the possibility to analyze the corresponding Bethe wave functions
and norm formulas as limits of Bethe wave type functions associated to 
arbitrary level. Such an approach has been taken in recent years to 
successfully analyze 
the Bethe vectors for Gaudin models and their norms 
as critical level limit of integral solutions to 
Knizhnik-Zamolodchikov equations, see \cite{MV} and references therein, which
in turn has interesting connections to the geometric Langlands
correspondence \cite{FFT}. 

Another important background for the current paper is the study of 
Knizhnik-Za\-mo\-lod\-chi\-kov-Bernard (KZB) heat equations, cf., e.g., 
cf. \cite{Be,EK1,FV1}. 
Correlation functions on the torus satisfy, besides KZB  
equations, the KZB heat equation.
For one-point correlation functions the KZB heat equation is the only
equation that remains. In that case
it is natural to view the KZB heat equation as a 
nonstationary Schr{\"o}dinger equation involving the 
modular parameter as time variable. 
At critical level it reduces to the quantum Hamiltonian 
of a vector valued
version of the quantum elliptic Calogero-Moser-Sutherland model
(see \cite[\S 4]{EK1} and \cite{FV1}). 
Special solutions are the so called affine Jack polynomials; 
they can be defined algebraically \cite[\S 6]{EK2} or be expressed 
\cite[Thm. 7.6]{EK2} in terms of one-point correlation functions.
Cherednik \cite[Thm. 4.3]{Cherednik1995} has derived
such nonstationary Schr{\"o}dinger equations
using the action of the trigonometric 
Cherednik algebra by infinite trigonometric Dunkl operators.
$q$-Extensions of various of these results have been obtained, see, e.g.,
\cite[\S 11]{EK2}, \cite{FV2,FV3} and \cite{Ch}.

Let us describe now the results presented in this paper. 
In \cite{EOS0,EOS} a suitable realization of
the trigonometric Cherednik algebra at critical level 
in terms of affine Dunkl operators involving Heaviside functions
enters the study of
vector valued root system analogs of the quantum Bose gas on the circle
with pairwise delta-function interactions
(Dunkl type operators involving Heaviside functions appeared before
in \cite{Po,MW,KH}). Following the idea of Cherednik \cite{Cherednik1995},
we generalize this realization to an extended version of the 
trigonometric Cherednik algebra at arbitrary level. It contains a
quadratic Casimir element which we show to produce an explicit 
nonstationary Schr{\"o}dinger equation involving delta-potentials
(see Proposition \ref{Hamiltonian}).

We construct solutions 
of this nonstationary Schr{\"o}dinger equation 
using an affine version of the coordinate Bethe ansatz. 
The associated generalized Bethe wave functions are defined
in terms of an explicit
cocycle of the extended affine Weyl group, 
which is obtained from the normalized intertwiners of 
the trigonometric Cherednik algebra.

We show that the generalized
Bethe wave function satisfies a consistent system of equations as function
of the momenta. These equations are expected to be 
degenerations of Cherednik's 
\cite{Ch} affine difference-elliptic quantum
Knizhnik-Zamolodchikov equations. These equations replace
the requirement for the Bethe wave function of the
(vector valued) root system version of the quantum Bose gas on the circle
with pairwise delta-function interactions that the momenta satisfy
Bethe ansatz type equations 
(see \cite[Thm. 2.6]{EOS0} and \cite[Thm. 5.10]{EOS}).

With suitable restrictions on the momenta
(see \eqref{convergencedomain} for the explicit requirements)
we show that the Bethe wave functions for the vector valued root system
analog of the quantum Bose gas on the line (see \cite{HO})
are limits of the generalized
Bethe wave functions. The limit to critical level is more
subtle. We will only make some preliminary 
comments on it in this paper. A thorough understanding of this limit
is expected to lead to new insights on the
root system analogs of the quantum Bose gas on the 
circle with pairwise delta-function potentials,
for instance on the quadratic
norms of the scalar Bethe wave functions. 
\vspace{.5cm}\\
\noindent
{\bf Acknowledgments}\\
The authors were partially supported by the Netherlands
Organization for Scientific Research (NWO) in the VIDI-project
``Symmetry and modularity in exactly solvable models''.
\section{The extended trigonometric Cherednik algebra}
\subsection{Notations}
Let $\gf$ be a complex finite dimensional reductive Lie algebra with
$\gf_s:=[\gf,\gf]$ simple. 
Write $\hf=\hf_s\oplus Z(\gf)$
with $\hf_s$ a Cartan subalgebra of $\gf_s$ and with
$Z(\gf)$ the center of $\gf$.

Let $(\cdot,\cdot)_s:\hf_s\times\hf_s\to\C$
be the restriction of the Killing form of $\gf_s$ to $\hf_s\times\hf_s$.
Let $R=R(\gf_s,\hf_s)\subset\hf_s$ be the set of roots of $\gf_s$ with respect
to $\hf_s$. It is a finite, reduced, irreducible crystallographic root system,
with the ambient Euclidean space $V_s$ taken to be the real span of $R$
and with scalar product the restriction of the bilinear form $(\cdot,\cdot)_s$
to $V_s\times V_s$. The root lattice $Q$ (respectively the co-root lattice
$Q^\vee$) is the rational integral span of all the roots $\al\in R$ 
(respectively the co-roots $\al^\vee$ ($\alpha\in R$)). 
These are lattices in $\hf$ satisfying 
$\mathbb{C}\otimes_{\mathbb{Z}}Q^\vee=\hf_s=\mathbb{C}\otimes_{\mathbb{Z}}Q$.

We fix a real form $V$ of $\hf$ of the form
$V_s\oplus V^\prime$ with
$V^\prime$ a real form of $Z(\gf)$. We 
extend the scalar product $(\cdot,\cdot)_s$ on $V_s$
to a scalar product $(\cdot,\cdot)$ on $V$
such that $V^\prime\perp V_s$. Its complex bilinear extension
to a 
bilinear form on $\hf$ is also denoted by $(\cdot,\cdot)$. 
We use it to identify the linear dual
$\hf^*$ with $\hf$. Set $\textup{O}(\hf)$ for the group 
of invertible complex linear endomorphisms of $\hf$ preserving the bilinear
form $(\cdot,\cdot)$.

Put $\widetilde{\hf}=\hf\oplus \C c$, 
$\wh\hf=\hf\oplus \C c\oplus\C d$, and extend the
form $(\cdot,\cdot)$ to a 
non-degenerate symmetric bilinear form on $\wh\hf$
by requiring 
\[
(c,d)=1,\quad (c,c)=(d,d)=(c,v)=(d,v)=0\qquad\forall\, v\in\hf.
\]

Write $\textup{O}(\wh\hf)\subset\GL(\wh\hf)$ for the subgroup of
invertible linear endomorphisms preserving the bilinear form $(\cdot,\cdot)$
on $\wh{\hf}$.
We identify $\textup{O}(\hf)$ with the subgroup
\[\{\sigma\in\textup{O}(\wh\hf)\,\, | \,\, \sigma(\hf)=\hf\,\,\,
\&\,\,\, \sigma(c)=c,\, \sigma(d)=d\}
\]
of $\textup{O}(\wh\hf)$.
For $u\in\hf$ define $t_{u}\in\textup{O}(\wh\hf)$
by
\begin{equation}\label{eq:Waaction_ty}
t_{u}(v+\eta c+\xi d)=v+\xi u + 
\big(\eta-\frac{\xi}{2}(u,u) -(v,u)\big)c+\xi d,
\end{equation}
where $v\in\hf$ and $\eta,\xi\in\mathbb{C}$. The map 
$\hf\rightarrow \textup{O}(\wh\hf)$ given by $u\mapsto t_u$ is a monomorphism
of groups. 
If $\sigma\in\textup{O}(\hf)$ then 
$\sigma\circ t_u=t_{\sigma(u)}\circ\sigma$. We conclude that
$\textup{O}(\wh\hf)$ naturally contains the subgroup 
$\textup{O}(\hf)\ltimes\hf$ of affine linear transformations of $\hf$.

For $a\in\wh\hf$ such that $(a,a)\not=0$ we
set $a^\vee:=2a/(a,a)\in\wh\hf$ and
\[
s_a(\widehat{v}):=\widehat{v}-(\widehat{v},a^\vee)a.
\]
Note that $s_a\in\textup{O}(\wh\hf)$ is an involution.
It fixes $c$ if $a\in\widetilde{\hf}$
and also $d$ if $a\in\hf$.
Write $W\subset \textup{O}(\hf)$ for the Weyl group of $R$ 
generated by $s_\al$ ($\al\in R$). The lattices $Q$ and $Q^\vee$ are
$W$-invariant. 

Set $\Ra=R\oplus\Z c\subseteq\widetilde{\hf}$. It plays the role
of the set of real roots of the 
untwisted affine Lie algebra associated to $\gf_s$, cf. \cite{Kac1990}.
If we identify
$\widetilde{V}=V\oplus\mathbb{R}c$ 
with the linear functionals on $V$ by
interpreting $v+\xi c$ ($v\in V$, $\xi\in\mathbb{R}$) as the affine linear
functional mapping $v^\prime$ to $(v,v^\prime)+\xi$, then $\Ra$ is a 
reduced irreducible affine root system in the sense of Macdonald
\cite{Mac1972}. 

The affine Weyl group $\Wa$ of $\Ra$ is defined to be
the subgroup of $\textup{O}(\wh\hf)$
generated by $s_a$ ($a\in\Ra$). The affine Weyl group is contained
in the subgroup of affine linear transformations of $\hf$ since
\[
s_{\alpha+mc}=s_\alpha t_{m\alpha^\vee},\qquad \alpha\in R,\, 
m\in\mathbb{Z}.
\]
In particular $\Wa\simeq W\ltimes Q^\vee$.

Fix a choice of positive roots $R^+$ of $R$ and let 
$F=\{a_1,\ldots,a_n\}$ be the associated set of simple roots. 
Let $\theta$ be the highest root in $R^+$ and set $a_0:=-\theta+c$.
Then $\widehat{F}:=\{a_0\}\cup F$ is a set of simple roots of $\Ra$,
i.e. every affine root $a\in\Ra$ can uniquely be written as
a nonnegative or nonpositive rational integral 
combination of the $a_i$'s. Denote $\Ra^{+}$ and $\Ra^-$
for the associated sets
of positive and negative affine roots, respectively. 
The affine Weyl group $\Wa$ is a Coxeter
group with Coxeter generators the simple reflections $s_i:=s_{a_i}$ 
($0\leq i\leq n$). The finite Weyl group $W\subset\Wa$ is the standard
parabolic subgroup generated by $s_i$ ($1\leq i\leq n$).
Note that $s_0=s_\theta t_{-\theta^\vee}=t_{\theta^\vee}s_\theta$. 

We fix a lattice $Y$ in $V$ containing $Q^\vee$ and satisfying 
$(Y,Q)\subseteq\Z$. Note that $Y$ is automatically $W$-invariant.
The associated extended affine Weyl group is defined 
by $\Wa^Y:=W\ltimes Y$. It contains the affine Weyl group $\Wa$ 
as a normal subgroup, and $\Wa^Y/\Wa\simeq Y/Q^\vee$ is abelian.

The affine root system $\Ra\subset\widetilde{\hf}$ and
the level $\xi$ hyperplanes 
\begin{equation}\label{hxi}
\wh\hf_\xi:=\widetilde{\hf}+\xi d
\end{equation} 
are $\Wa^Y$-invariant. Furthermore, $\Wa^Y\simeq\Omega\ltimes\Wa$ with 
$\Om=\Om^Y$ the subgroup of $\Wa^Y$ consisting of elements $\om$
such that $\om(\Rap)\subseteq\Rap$. The group $\Om$ permutes the simple
affine roots.
\subsection{The definition}

The trigonometric Cherednik algebra, also known as the
degenerate double affine Hecke algebra, was defined in
\cite[Def. 1.1]{Cherednik1995}. We use an extended version, defined as follows.
\begin{defi}
Let $k:\Ra\to\mathbb{C}$ be a $\Wa^Y$-invariant function, which
we call a multiplicity function.
The extended trigonometric Cherednik algebra $\wh{H}^Y(k)=\wh{H}^Y(k;V)$
is the associative unital $\C$-algebra satisfying
\begin{enumerate}
\item $\wh{H}^Y(k)$ contains the symmetric algebra $S(\wh\hf)$
and the group algebra $\C[\Wa^Y]$ as subalgebras,
\item the multiplication map defines a linear isomorphism 
$S(\wh\hf)\otimes_\C \C[\Wa^Y]\simeq \wh{H}^Y(k)$,
\item the following cross relations hold:
\begin{equation}\label{cr}
\begin{split}
s_a\cdot \wh v&=s_a(\wh v)\cdot s_a - k_a (a,\wh v), 
\qquad\forall a\in \Fa,\; \forall\,\wh v\in \wh\hf,\\
\om\cdot \wh v&=\om(\wh v)\cdot \om, 
\qquad\qquad\qquad\,\,\,\forall\om\in\Om,\; \forall\, \wh v\in \wh\hf.
\end{split}
\end{equation}
\end{enumerate}
The subalgebra $\widetilde{H}^{Y}(k)$ of $\wh{H}^Y(k)$ generated by
$S(\widetilde{\hf})$ and $\C[\Wa^Y]$ is the
trigonometric Cherednik algebra.
\end{defi}
\begin{rema}\label{notfull}
{\bf (i)} In the next subsection we verify that
$\wh{H}^Y(k)$ is well defined (see Proposition \ref{welldefined}).\\
{\bf (ii)}
The trigonometric Cherednik algebra $\widetilde{H}^{Y}(k)$ admits a similar
characterization (1)-(3) 
as $\wh{H}^Y(k)$, with the role of the symmetric algebra
$S(\wh\hf)$ taken over by $S(\widetilde{\hf})$.\\
{\bf (iii)} 
Let $U\subseteq V$ be the real span of $Y$ and $U^\perp_{\mathbb{C}}
\subset\hf$ the complexification of the orthocomplement $U^\perp$
of $U$ in $V$. Then
$\wh{H}^Y(k;V)\simeq\wh{H}^Y(k;U)\otimes_{\mathbb{C}}S(U^\perp_{\mathbb{C}})$ as
algebras.
\end{rema}
We write $\wh{H}(k):=\wh{H}^{Q^\vee}(k)$ and $\widetilde{H}(k):=
\widetilde{H}^{Q^\vee}(k)$. Note that
\[
\wh{H}^Y(k)\simeq \Omega\ltimes\wh{H}(k),\qquad
\widetilde{H}^{Y}(k)\simeq\Omega\ltimes\widetilde{H}(k).
\]
Observe that $c\in Z(\wh{H}^Y(k))$. 
We define the \emph{extended trigonometric Cherednik algebra at level  
$\kappa$} to be $\wh{H}_{\kappa}^Y(k) := \wh{H}^Y(k)/(c-\kappa)$.

We call $\widetilde{H}_{\kappa}^{Y}(k)
:=\widetilde{H}^{Y}(k)/(c-\kappa)$ the trigonometric Cherednik algebra
at level $\kappa$, and $\widetilde{H}_0^Y(k)$ the 
trigonometric Cherednik algebra at critical level. 
A detailed analysis of $\widetilde{H}_0^{Y}(k)$ in connection to the
root system analogs of the one dimensional quantum Bose gas with pairwise
delta-function interactions was undertaken 
in \cite{EOS0} and \cite{EOS}. In this context the quantum Hamiltonians
of the quantum integrable system arise from the center 
$Z(\widetilde{H}_0^{Y}(k))$ of $\widetilde{H}_0^{Y}(k)$, which 
contains $S(\hf)^W$. In this paper we consider these structures
away from critical level, replacing the role of $\widetilde{H}_0^{Y}(k)$
by the extended trigonometric Cherednik algebra $\wh{H}^Y(k)$.
Its center is described as follows. 
\begin{prop}\label{centerprop}
If $Y$ is a full lattice in $V$ then
$Z\bigl(\wh{H}^Y(k)\bigr)=S(\wh{\hf})^{\Wa^Y}=
\mathbb{C}[c,C]$, where $C\in S(\wh{\hf})^{\Wa^Y}$ is the
Casimir element
\begin{equation}\label{Casimir}
C=\sum_{i=1}^mv_i^2+2cd,
\end{equation}
and $\{v_i\}_{i=1}^m$ is an (arbitrary choice of) orthonormal basis
of $V$ with respect to $(\cdot,\cdot)$.
\end{prop}
\begin{proof}
Let $\{b_j\}_{j=1}^{m+2}$ be a basis of $\wh{V}$ and
$\{b^j\}_{j=1}^{m+2}$ its dual basis with respect to $\bigl(\cdot,\cdot\bigr)$.
The Casimir element can alternatively be defined by
\begin{equation}\label{Casimiralt}
C=\sum_{j=1}^{m+2}b_jb^j.
\end{equation}
{}From the $\Wa^Y$-invariance of $\bigl(\cdot,\cdot\bigr)$ it then follows
that $C\in S(\wh{\hf})^{\Wa^Y}$. 

Now take an orthonormal basis $\{v_i\}_{i=1}^n$ 
of $V_s$. Set $\wh{\hf}_s=\hf_s\oplus\mathbb{C}c
\oplus\mathbb{C}d$ and write
\[
C_s:=\sum_{i=1}^nv_i^2+2cd\in S(\wh{\hf}_s).
\]
By \cite[Prop. 4.1]{CI} we have 
$S(\wh{\hf}_s)^{\Wa}=\mathbb{C}[c,C_s]$. Considering
$S(\wh{\hf}_s)$ and $S(Z(\mathfrak{g}))$ as $\Wa$-module subalgebras of 
$S(\wh{\hf})$, we conclude that
\[
S(\wh{\hf})^{\Wa}=S(\wh{\hf}_s)^{\Wa}S(Z(\mathfrak{g}))=
\mathbb{C}[c,C_s]S(Z(\mathfrak{g}))=
\mathbb{C}[C]S(Z(\mathfrak{g})\oplus\mathbb{C}c),
\]
where the last equality follows from the fact that
$C-C_s\in S(Z(\mathfrak{g}))$. Set $Y_0:=Y\cap V^\prime$ (recall that 
$V^\prime$ is the real form of $Z(\mathfrak{g})$ such 
that $V=V_s\oplus V^\prime$).
Since $C$ is $\Wa^Y$-invariant we conclude that
\[
S(\wh{\hf})^{\Wa^Y}\subseteq 
\mathbb{C}[C]S(Z(\mathfrak{g})\oplus\mathbb{C}c)^{Y_0}.
\]
Furthermore,
\[
S(Z(\mathfrak{g})\oplus\mathbb{C}c)^{Y_0}=\mathbb{C}[c]
\]
since $Y_0$ is a full lattice in $V^\prime$,
cf. the proof of \cite[Prop. 4.1]{CI}. We conclude that
$S(\wh{\hf})^{\Wa^Y}=\mathbb{C}[c,C]$. 
Hence it remains to show that
\[
Z(\wh{H}^Y(k))=S(\wh{\hf})^{\Wa^Y}.
\]
This can be proved 
by a straightforward adjustment of 
the analogous statement (due to Lusztig)
for the degenerate affine Hecke algebra, cf. \cite[Prop. 1.1]{Op}. 
\end{proof}

\subsection{Difference-reflection operators}
We construct now a representation of the extended trigonometric Cherednik
algebra $\wh{H}^Y(k)$ using Dunkl operators involving Heaviside functions.
It generalizes the representation of the trigonometric
Cherednik algebra $\widetilde{H}_0^Y(k)$ at critical level
constructed in \cite[\S 4.2]{EOS}.

Set $\widetilde{V}=V\oplus \mathbb{R}c$ and
$\widehat{V}:=\widetilde{V}\oplus\mathbb{R}d$.
The bilinear form $(\cdot,\cdot)$ on $\wh\hf$ restricts to a real valued
non-degenerate symmetric bilinear form on $\widehat{V}$.
We also write $\widehat{V}^+:=\widetilde{V}\oplus\mathbb{R}_{>0}d$
and $\widehat{V}_\xi:=\widetilde{V}+\xi d$. Note that $\widetilde{V}$,
$\widehat{V}$, $\widehat{V}^+$ and $\widehat{V}_\xi$ are $\wh{W}^Y$-invariant
subsets of $\wh{\hf}$.

The open subset
\[
\widehat{V}_{reg}^+:=\{\widehat{v}=v+\eta c+\xi d\in\widehat{V}^+\,\, 
| \,\, (a,\widehat{v})\not=0\quad \forall
a\in\widehat{R}^+\}
\]  
of regular elements in $\wh{V}^+$ is $\wh{W}^Y$-invariant.
Denote by $\mathcal{C}$ the collection of connected
components of $\widehat{V}_{reg}^+$. The affine Weyl group $\Wa$ acts
simply transitively on $\mathcal{C}$.
The convex polytope
\begin{equation*}
\begin{split}
\widehat{C}_+:=&\{ \widehat{v}\in\widehat{V}^+\,\, | \,\, 
(a,\widehat{v})>0\quad \forall\, a\in\widehat{F} \}\\
=&\{\widehat{v}\in V_s\oplus\mathbb{R}_{>0}d\,\, | \,\, (a,\widehat{v})>0
\quad \forall\, a\in\widehat{F} \}\oplus V^\prime\oplus \mathbb{R}c
\end{split}
\end{equation*}
is a connected component of $\widehat{V}_{reg}^+$
which we call the fundamental chamber. 
Note that $\om(\widehat{C}_+)=
\widehat{C}_+$ for all $\omega\in\Omega$.

$M$ will always stand for a finite dimensional, 
complex, left $\Wa^Y$-module. Its representation map
will be denoted by $\pi_M$.

We define the complex vector space
\[
\mathcal{F}_M:=\prod_{\wh{C}\in\mathcal{C}}
\bigl(C^\omega(\widehat{V}^+)\otimes_{\mathbb{C}}M\bigr),
\]
where $C^\omega(\widehat{V}^+)$ is the space of complex valued, real analytic
functions on $\widehat{V}^+$. 
An element $f=(f_{\wh{C}})_{\wh{C}\in\mathcal{C}}\in\mathcal{F}_M$
should be thought of as a collection of real analytic 
$M$-valued functions $f_{\wh{C}}$ on $\wh{C}$ ($\wh{C}\in\mathcal{C}$)
with the additional requirement that each $f_{\wh{C}}$ admits 
a real analytic extension to $\widehat{V}^+$. We define
the support $\textup{supp}(f)$ of $f=(f_{\wh{C}})_{\wh{C}\in\mathcal{C}}
\in\mathcal{F}_M$ to
be the collection of connected components $\wh{C}$ for which $f_{\wh{C}}$ 
is nonzero. Note that $\Wa^Y$
acts on $f=(f_{\wh{C}})_{\wh{C}\in\mathcal{C}}\in\mathcal{F}_M$ by 
\begin{equation}\label{actionF}
(wf)_{\wh{C}}(\widehat{v}):=
\pi_M(w)\bigl(f_{w^{-1}\wh{C}}(w^{-1}\widehat{v})\bigr),
\qquad w\in\Wa^Y,\, \wh{C}\in\mathcal{C},\, \widehat{v}\in\widehat{V}^+.
\end{equation}
For $\widehat{v}\in\widehat{V}$ we define the linear endomorphism 
$\partial_{\widehat{v}}$ of $\mathcal{F}_M$ as the componentwise directional
derivative, 
\[
(\partial_{\widehat{v}}f)_{\wh{C}}(\wh{u}):=(\partial_{\widehat{v}}f_{\wh{C}})(\wh{u})
=\frac{d}{dt}|_{t=0}f_{\wh{C}}(\wh{u}+t\wh{v}).
\]
For $a\in\widehat{R}^+$ let $\chi_a: \widehat{V}\rightarrow \{0,1\}$
be the characteristic function of the half-space 
$H_a^-:=\{\widehat{v}\in\widehat{V}\,\, | \,\, (a,\widehat{v})<0\}$.
For each chamber $\wh{C}\in\mathcal{C}$ either $\chi_a|_{\wh{C}}\equiv 1$
or $\chi_a|_{\wh{C}}\equiv 0$.
We also write $\chi_a$ for the linear endomorphism
of $\mathcal{F}_M$ mapping $f=(f_{\wh{C}})_{\wh{C}\in\mathcal{C}}$ to
$\chi_af:=\{\chi_a(\wh{C})f_{\wh{C}}\}_{\wh{C}\in\mathcal{C}}$. 
\begin{lem}
Fix $\wh{C}=w(\wh{C}^+)\in\mathcal{C}$ ($w\in \wh{W}$). If $a\in\wh{R}^+$ then
the following two statements are equivalent:
\begin{enumerate}
\item $\wh{C}\in\textup{supp}(\chi_af)$ for some 
$f\in\mathcal{F}_M$, 
\item $a\in\wh{R}^+\cap w(\wh{R}^-)$ (which is a finite
set of positive affine roots).
\end{enumerate}
\end{lem}
\begin{proof}
Both are easily seen to be equivalent to $\chi_a(\wh{C})=1$.
\end{proof}
The lemma allows us to define
linear operators $\mathcal{D}_{\widehat{v}}^M$ ($\widehat{v}\in
\widehat{V}$) on $\mathcal{F}_M$ by 
\[
\mathcal{D}_{\widehat{v}}^Mf:=\partial_{\widehat{v}}f-
\sum_{a\in\Ra^+}k_a(a,\widehat{v})\chi_as_af.
\]

The following proposition extends the results from
\cite[\S 4.2]{EOS}. 
\begin{prop}\label{welldefined}
Let $k: \Ra\rightarrow\mathbb{C}$ be a multiplicity function.
\begin{enumerate}
\item The extended trigonometric Cherednik algebra $\wh{H}^Y(k)$
is well defined.
\item The assignments
\begin{equation*}
\begin{split}
\widehat{v}&\mapsto \mathcal{D}_{\widehat{v}}^M,\qquad 
\widehat{v}\in\widehat{V},\\
w&\mapsto w,\qquad\,\,\, w\in\Wa^Y,
\end{split}
\end{equation*}
uniquely define an algebra morphism $\widehat{\pi}: \wh{H}^Y(k)\rightarrow
\textup{End}_{\mathbb{C}}\bigl(\mathcal{F}_M\bigr)$. 
\end{enumerate}
\end{prop}
\begin{proof} 
Repeating the arguments of the proof of \cite[Thm. 4.1]{EOS0}
gives that the operators $\mathcal{D}_{\widehat{v}}^M$ 
($\widehat{v}\in\widehat{V}$), $s_a$ ($a\in\wh{F}$) and $\om$ ($\om\in\Omega$)
on $\mathcal{F}_M$ satisfy the defining relations
of $\wh{H}^Y(k)$,
\begin{equation*}
\begin{split}
s_a\mathcal{D}_{\widehat{v}}^M&=\mathcal{D}_{s_a\widehat{v}}^Ms_a-
k_a(a,\widehat{v}),\\
\omega\mathcal{D}_{\widehat{v}}^M&=\mathcal{D}_{\omega\widehat{v}}^M\,\omega,\\
\mathcal{D}_{\widehat{v}}^M\mathcal{D}_{\widehat{v}^\prime}^M&=
\mathcal{D}_{\widehat{v}^\prime}^M\mathcal{D}_{\widehat{v}}^M.
\end{split}
\end{equation*}
For $Y=Q^\vee$ and $M=\textup{Triv}$ the trivial one-dimensional $\Wa$-module
the resulting linear map
$S(\wh\hf)\otimes_{\mathbb{C}}
\mathbb{C}[\Wa]\rightarrow\textup{End}(\mathcal{F}_{\textup{Triv}})$
is easily seen to be injective,
hence the extended trigonometric Cherednik algebras $\wh{H}(k)$ and
$\wh{H}^Y(k)=\Omega\ltimes\wh{H}(k)$ are well defined. Part
(2) of the proposition follows now immediately.
\end{proof}
\begin{rema}
A trigonometric version of the representation
$\wh{\pi}|_{\widetilde{H}^Y(k)}$ was constructed 
by Cherednik in \cite[Thm. 3.1]{Cherednik1995} using
infinite trigonometric
Dunkl operators. He considers separately
a trigonometric analog of the operator $\wh{\pi}(d)$ 
(see \cite[(4.12)]{Cherednik1995}). He
remarks that the cross relations \eqref{cr} are respected
(see \cite[(4.13)]{Cherednik1995}) but that, in contrast to our
setup, it does not result in a
representation of the extended trigonometric Cherednik algebra $\wh{H}^Y(k)$.
\end{rema}
Note that $\mathcal{D}_c^M=\partial_c$ since $(a,c)=0$ for all
$a\in\wh{R}$. Hence for $\kappa\in\mathbb{C}$,
\[
\mathcal{F}_M^\kappa:=\{f\in\mathcal{F}_M \,\, | \,\, \partial_cf=\kappa f\}
\]
is a $\wh{H}^Y(k)$-submodule of $\mathcal{F}_M$, and the action of
$\wh{H}^Y(k)$ on $\mathcal{F}_M^\kappa$ descends to an
action of the extended trigonometric Cherednik algebra
$\wh{H}^Y_\kappa(k)$ at level $\kappa$. 

For $\wh{\lambda}\in\wh{\hf}$
let $e^{\wh{\lambda}}\in\mathcal{O}(\wh{\hf})$
be the holomorphic function $\wh{\mu}\mapsto
e^{(\wh{\lambda},\wh{\mu})}$. Its restriction to $\wh{V}^+$ defines
a complex valued, real analytic function on $\wh{V}^+$ which we
also will denote by $e^{\wh{\lambda}}$.
We have
\[\mathcal{F}_M^\kappa=\prod_{\wh{C}\in\mathcal{C}}\bigl(
C^\omega_\kappa(\wh{V}^+)\otimes_{\mathbb{C}}M\bigr)
\]
where
\[
C_\kappa^\omega(\widehat{V}^+):=
\{f\in C^\omega(\widehat{V}^+)\,\, | \,\, \partial_c(f)=\kappa f\}.
\]
Note that $C_\kappa^\omega(\widehat{V}^+)=
e^{\kappa d}C_0^\omega(\widehat{V}^+)$.

\subsection{Integral-reflection operators}
In this subsection we 
give a second representation of the extended trigonometric Cherednik
algebra $\wh{H}^Y(k)$, now in terms of integral-reflection operators.
The results in this subsection build on constructions from 
\cite{Gutkin1982,HO,EOS0,EOS}.

For $a\in\widehat{R}$ define an integral-reflection operator
$I(a)$ on $C^\omega(\widehat{V}^+)$ by
\begin{equation}\label{integralreflection}
(I(a)f)(\widehat{v}):=\int_0^{(a,\widehat{v})}f(\widehat{v}-ta^\vee)dt.
\end{equation}
The following is a version of \cite[Thm. 4.11]{EOS}
at unspecified level. 
\begin{thm}\label{irthm}
Let $M$ be a finite dimensional $\Wa^Y$-module.
There exists a unique algebra homomorphism $\widehat{Q}: 
\wh{H}^Y(k)\rightarrow \textup{End}_{\mathbb{C}}\bigl(C^\omega(\widehat{V}^+)
\otimes_{\mathbb{C}}M\bigr)$ satisfying
\begin{equation}\label{explQ}
\begin{split}
\widehat{Q}(s_a)&=s_a\otimes\pi_M(s_a)-k_aI(a)\otimes\textup{Id}_M,\qquad 
\, a\in\widehat{F},\\
\widehat{Q}(\omega)&=\omega\otimes\pi_M(\omega),\quad\qquad\qquad\qquad\qquad
\,\,\,\omega\in\Omega,\\
\widehat{Q}(\widehat{v})&=\partial_{\widehat{v}}\otimes\textup{Id}_M,
\,\,\qquad\qquad\qquad\qquad\qquad \widehat{v}\in\widehat{V}.
\end{split}
\end{equation}
\end{thm}
\begin{proof}
Consider $S(\wh{\hf})\otimes_{\mathbb{C}}M^*$ as left $\wh{H}^Y(k)$-module
by the canonical vector space identification
\[
\textup{Ind}_{\mathbb{C}[\wh{W}^Y]}^{\wh{H}^Y(k)}(M^*)\simeq
S(\wh{\hf})\otimes_{\mathbb{C}}M^*.
\]
Using the complex linear antiinvolution ${}^\dagger: \wh{H}^Y(k)
\overset{\sim}{\longrightarrow}\wh{H}^Y(k)$ defined by $w^\dagger=w^{-1}$
($w\in\wh{W}^Y$) and $\wh{v}^\dagger=\wh{v}$ ($\wh{v}\in\wh{V}$),
its linear dual $(S(\wh{\mathfrak{h}})\otimes_{\mathbb{C}}M^*)^*$ becomes
a left $\wh{H}^Y(k)$-module. 

View 
$C^\omega(\wh{V}^+)\otimes_{\mathbb{C}}M$ as linear subspace of 
$(S(\wh{\hf})\otimes_{\mathbb{C}}M^*)^*$ by interpreting
$f\otimes m$ ($f\in C^\omega(\wh{V}^*)$,
$m\in M$) as the
linear functional 
\[
p\otimes\psi\mapsto \psi(m)(p(\partial)f)(0),\qquad
p\in S(\wh{\hf}),\, \psi\in M^*
\]
on $S(\wh{\hf})\otimes_{\mathbb{C}}M^*$,
where $p(\partial)$ is the constant coefficient partial differential
operator associated to $p\in S(\wh{\hf})$. 
Then 
$C^\omega(\wh{V}^+)\otimes_{\mathbb{C}}M$ is a $\wh{H}^Y(k)$-submodule
of $(S(\wh{\hf})\otimes M^*)^*$. A direct computation establishes
the explicit formulas \eqref{explQ} for the resulting action of
$\wh{H}^Y(k)$ on $C^\omega(\wh{V}^+)\otimes_{\mathbb{C}}M$.
\end{proof}

\begin{rema}
The proof of Theorem \ref{irthm} is simpler than the proof
at critical level (see \cite[Thm. 4.11]{EOS}), since at critical level
the arguments of the above proof lead to an explicit action
by integral-reflection operators 
which is not yet of the desired form (see \cite[Cor. 4.10]{EOS}).
\end{rema}

The integral-reflection operators $\widehat{Q}(w)$ ($w\in\Wa$)
gives rise to a linear map
\[
T: C^\omega(\widehat{V}^+)\otimes_{\mathbb{C}}M\rightarrow\mathcal{F}_M,
\]
with $Tf=\{(Tf)_C\}_{C\in\mathcal{C}}\in\mathcal{F}_M$ 
for $f\in C^\omega(\widehat{V}^+)\otimes_{\mathbb{C}}M$
defined by
\[
(Tf)_{w\widehat{C}^+}(\widehat{v}):=
\pi_M(w)\bigl(\bigl(\widehat{Q}(w^{-1})f\bigr)(w^{-1}\widehat{v})\bigr),
\qquad w\in\Wa,\,\,\widehat{v}\in\widehat{V}^+.
\]
It is the unique linear map satisfying $(Tf)_{\widehat{C}^+}=f$
and $T\circ\widehat{Q}(w)=\widehat{\pi}(w)\circ T$ for all $w\in\Wa$, cf. 
\cite[Lemma 4.14]{EOS}. In fact, in analogy with
\cite[Prop. 4.15]{EOS} we have
\begin{prop}\label{propagationproperty}
The linear map $T: C^\omega(\widehat{V}^+)\otimes_{\mathbb{C}}M\rightarrow
\mathcal{F}_M$ is $\wh{H}^Y(k)$-linear,
\[
T\circ\wh{Q}(h)=\wh{\pi}(h)\circ T,\qquad \forall\, h\in\wh{H}^Y(k).
\]
\end{prop}
\begin{rema}
For $\kappa\in\mathbb{C}$ the space $C_\kappa^\omega(\wh{V}^+)\otimes M$
is a $\wh{H}^Y(k)$-submodule of $C^\omega(\wh{V}^+)\otimes_{\mathbb{C}}M$
with respect to the $\wh{Q}$-action. The action on
$C_\kappa^\omega(\wh{V}^+)\otimes M$
descends to an action of $\wh{H}^Y_\kappa(k)$
since $\wh{Q}(c)=\partial_c$.
The intertwiner $T$ restricts to an intertwiner
\[
T_\kappa: C^\omega_\kappa(\wh{V}^+)\otimes_{\mathbb{C}}M\rightarrow
\mathcal{F}_M^\kappa
\]
of $\wh{H}^Y_\kappa(k)$-modules.
\end{rema}

\section{Generalized Bethe wave functions}

In this section we construct generalized Bethe wave functions, being
$\wh{Q}(\Wa^Y)$-invariant 
eigenfunctions of the constant coefficient differential operator $\wh{Q}(C)$,
as infinite series expansions of plane waves with explicit
cocycle values as coefficients. 
The relevant cocycle of $\Wa^Y$
arises from the normalized intertwiners of the trigonometric Cherednik
algebra $\widetilde{H}^Y(k)$. On the other hand, this 
cocycle can be used to define a 
consistent system of equations. These equations should be thought of as
formal degenerations
of Cherednik's \cite{Ch} affine difference-elliptic quantum
Knizhnik-Zamolodchikov equations. We show that the generalized Bethe
wave functions satisfy this consistent system of equations
as function of the momenta.

\subsection{The cocycle}
Let $k: \Ra\rightarrow\mathbb{C}$ be a $\Wa^Y$-invariant multiplicity function.
Set
\begin{equation*}
\begin{split}
\wh{\hf}_{reg}:=&\{\widehat{\lambda}\in\wh\hf \,\, | \,\,
(a^\vee,\widehat{\lambda})\not=0,k_a\quad \forall a\in\Ra\},\\
\wh{\hf}_{\kappa,reg}:=&\wh{\hf}_{reg}\cap\wh{\hf}_{\kappa},
\end{split}
\end{equation*}
so that $\wh{\hf}_{\kappa,reg}=\hf_{\kappa,\reg}+\mathbb{C}c+\kappa d$
with 
\[
\hf_{\kappa,reg}=\{\lambda\in\hf\,\, | \,\, (\alpha^\vee,\lambda)+
\frac{2m\kappa}{(\alpha,\alpha)}\not=0,
k_{\alpha+m c}\quad \forall\,\alpha\in R,\,\forall\, m\in\mathbb{Z}\}.
\]
Let $\mathbb{C}_+=\{z\in\mathbb{C}\,\, | \,\, \textup{Re}(z)>0\}$ be the open 
right half plane in $\mathbb{C}$.
Set $\wh{\hf}^+:=\widetilde{\hf}+\mathbb{C}_+d$ and 
\[
\wh{\hf}_{reg}^+:=\wh{\hf}_{reg}\cap\wh{\hf}^+=
\bigcup_{\kappa\in\mathbb{C}_+}\wh{\hf}_{\kappa,reg}.
\]
Then 
$\wh{\hf}_{reg}^+\subset\wh{\hf}$ is open, connected, and
the boundary of $\wh{\hf}_{reg}^+$ contains $\wh{\hf}_{0,reg}$.

Let $\mathbb{C}[\wh{\hf}_{reg}]$ be the subalgebra of the field of
rational functions on $\wh{\hf}$ obtained by adjoining $(a^\vee-k_a)^{-1}$
and $a^{-1}$ to $\mathbb{C}[\wh{\hf}]$ for all $a\in\Ra$. 
By \cite{Cherednik1997} we have
\begin{prop}
There exist unique $J_w\in\mathbb{C}[\wh{\hf}_{reg}]\otimes\mathbb{C}[\Wa^Y]$
($w\in\Wa^Y$) satisfying, as rational $\mathbb{C}[\Wa^Y]$-valued
functions on $\wh{\hf}_{reg}$, 
\[
J_{uw}(\widehat{\lambda})=J_u(w\widehat{\lambda})J_w(\widehat{\lambda})
\qquad \forall\, u,w\in\Wa^Y
\]
and satisfying
\begin{equation*}
\begin{split}
J_{s_a}(\widehat{\lambda})&=\frac{(a^\vee,\widehat{\lambda})s_a+k_a}
{(a^\vee,\widehat{\lambda})-k_a},\qquad a\in\widehat{F},\\
J_{\omega}(\widehat{\lambda})&=\omega,\qquad\qquad\qquad\quad\,\, 
\omega\in\Omega.
\end{split}
\end{equation*}
\end{prop}
The proof of the proposition
uses the normalized intertwiners of the trigonometric Cherednik algebra,
cf. \cite{Cherednik1997}.
\begin{rema}
Viewing $J_w$ as $\mathbb{C}[\Wa^Y]$-valued rational function on
$\wh{\hf}$, we have $\partial_c(J_w)=0$ for all $w\in\Wa^Y$.
\end{rema}
The ring $\mathcal{O}(\wh{\hf}_{reg}^+)$ of holomorphic functions
on $\wh{\hf}_{reg}^+$
is naturally a $\mathbb{C}[\wh{\hf}_{reg}]$-module.
If $M$ is a $\Wa^Y$-module, the cocycle $\{J_w\}_{w\in\Wa^Y}$ thus canonically
acts on $\mathcal{O}(\wh{\hf}_{reg}^+)\otimes_{\mathbb{C}}
\textup{End}_{\mathbb{C}}(M)$. 
\begin{cor}\label{action}
Let $M$ be a finite dimensional $\Wa^Y$-module.
Set
\begin{equation}\label{Psiw}
(\Psi\cdot w)(\widehat{\lambda}):=
\Psi(w\widehat{\lambda})J_{w}(\widehat{\lambda}),\qquad
w\in\Wa^Y,\, \Psi\in\mathcal{O}(\wh{\hf}_{reg}^+)
\otimes_{\mathbb{C}}\textup{End}_{\mathbb{C}}(M),
\end{equation}
where $\Psi$ is viewed as 
$\textup{End}_{\mathbb{C}}(M)$-valued holomorphic function on
$\wh{\hf}_{reg}^+$. Then \eqref{Psiw} defines a right $\Wa^Y$-action on
$\mathcal{O}(\wh{\hf}_{reg}^+)\otimes_{\mathbb{C}}\textup{End}_{\mathbb{C}}(M)$.
\end{cor}
We call the set of equations
\begin{equation}\label{dKZ}
\Psi(t_y\widehat{\lambda})J_{t_y}(\widehat{\lambda})=\Psi(\widehat{\lambda})
\qquad \forall\, y\in Y
\end{equation}
the affine difference Knizhnik-Zamolodchikov (adKZ) equations and 
\[
\textup{adKZ}:=\bigl(\mathcal{O}(\wh{\hf}_{reg}^+)\otimes_{\mathbb{C}}
\textup{End}_{\mathbb{C}}(M)\bigr)^{\cdot Y}
\]
the corresponding space of solutions. It
is a $\cdot W$-submodule of
$\textup{O}(\wh{\hf}_{reg}^+)\otimes_{\mathbb{C}}\textup{End}_{\mathbb{C}}(M)$.
\begin{rema}
The affine difference KZ equations \eqref{dKZ} are expected to be formal
degenerations of Cherednik's \cite{Ch} affine difference-elliptic 
quantum affine KZ equations, which are naturally associated to double
affine Hecke algebras.
\end{rema}
Explicitly, the affine difference KZ equations \eqref{dKZ} read
\begin{equation}\label{dKZexplicit}
\Psi\bigl(\lambda+\kappa y+(\eta-\frac{\kappa}{2}(y,y)-(\lambda,y))c+\kappa d)
J_{t_y}(\lambda+\kappa d)=\Psi(\lambda+\eta c+\kappa d)\qquad
\forall\, y\in Y.
\end{equation}
Set $\mathcal{O}_\xi(\wh{\hf}^+_{reg})=
\{f\in\mathcal{O}(\wh{\hf}^+_{reg})\,\, | \,\,
\partial_c(f)=\xi f\}$, so that $\mathcal{O}_\xi(\wh{\hf}_{reg}^+)=
e^{\xi d}\mathcal{O}_0(\wh{\hf}_{reg}^+)$.
Observe that 
$\mathcal{O}_\xi(\wh{\hf}_{reg}^+)\otimes_{\mathbb{C}} 
\textup{End}_{\mathbb{C}}(M)$ is a $\cdot\Wa^Y$-submodule
of $\mathcal{O}(\wh{\hf}_{reg}^+)\otimes_{\mathbb{C}}
\textup{End}_{\mathbb{C}}(M)$. Set
\[
\textup{adKZ}_\xi:=\bigl(\mathcal{O}_\xi(\wh{\hf}_{reg}^+)\otimes_{\mathbb{C}}
\textup{End}_{\mathbb{C}}(M)\bigr)^{\cdot Y}.
\]
By \eqref{dKZexplicit}, if $\Psi\in\mathcal{O}_\xi(\wh{\hf}_{reg}^+)
\otimes_{\mathbb{C}}\textup{End}_{\mathbb{C}}(M)$ then 
$\Psi\in\textup{adKZ}_\xi$ iff
\begin{equation}\label{dKZexplicitxi}
\Psi(\lambda+\kappa y+\eta c+\kappa d)e^{-\frac{\kappa\xi}{2}(y,y)-
\xi(\lambda,y)}J_{t_y}(\lambda+\kappa d)=\Psi(\lambda+\eta c+\kappa d)\qquad
\forall\, y\in Y
\end{equation}
as $\textup{End}_{\mathbb{C}}(M)$-valued holomorphic function in
$\lambda+\eta c+\kappa d\in\wh{\hf}_{reg}^+$. 
Note that \eqref{dKZexplicitxi} formally makes sense if $\kappa=0$,
in which case it gives
\begin{equation}\label{dKZexplicitxikappazero}
\Phi(\lambda)e^{-\xi(\lambda,y)}J_{t_y}(\lambda)=\Phi(\lambda),\qquad
\forall y\in Y
\end{equation}
for $\Phi(\lambda)\in\textup{End}_{\mathbb{C}}(M)$ 
and $\lambda\in\wh{\hf}_{0,reg}$. These equations are closely related to 
the Bethe ansatz equations for the vector valued
root system analogs of the quantum Bose gas on the circle with 
pairwise delta-function
interactions, see \cite[Thm. 5.10]{EOS}. 

For $\xi=0$ the affine difference KZ equations \eqref{dKZexplicitxi}
for fixed $\kappa\in\mathbb{C}_+$ take the form
\begin{equation}\label{dKZexplicitxizero}
\Phi(\lambda+\kappa y)J_{t_y}(\lambda+\kappa d)=\Phi(\lambda),\qquad
\forall\, y\in Y
\end{equation}
for $\Phi$ a $\textup{End}_{\mathbb{C}}(M)$-valued holomorphic function on
$\hf_{\kappa,reg}$. These are 
degenerations of Cherednik's \cite{Cherednik1992} quantum affine  
KZ equations associated to double
affine Hecke algebras. They form a consistent system of 
difference equations naturally compatible to 
trigonometric KZ equations, see \cite{Stokman2011}.
\subsection{Bethe wave functions at critical level}
Before constructing generalized Bethe wave functions and their relation
to the affine difference KZ equations,
we first recall the related results at critical level from \cite{EOS}.
Let $\xi>0$ and write
$Q_\xi: \widetilde{H}_0^{Y}(k)\rightarrow\textup{End}_{\mathbb{C}}(C^\omega(V)
\otimes_{\mathbb{C}}M)$ for the analog of the integral-reflection action
$\wh{Q}$ (see \cite[Thm. 4.11]{EOS}), defined by
$Q_\xi(v)=\partial_v$ ($v\in V$) and 
\begin{equation*}
\begin{split}
(Q_\xi(s_a)f)(v)&=
\pi_M(s_a)f(s_a\circ_\xi v)-k_a\int_0^{a^\xi(v)}f(v-tDa^\vee)dt,\qquad
\,\,\, a\in\wh{F},\\
(Q_\xi(\omega)f)(v)&=\pi_M(\omega)f(\omega^{-1}\circ_\xi v),\,\,\,\,
\quad\qquad\qquad
\qquad\qquad\quad\quad\qquad \omega\in\Omega,
\end{split}
\end{equation*}
where $a^\xi(v)=(\alpha,v)+m\xi$ and $Da=\alpha$ for $a=\alpha+mc$ 
and $v\in V$, and where the action $\wh{W}^Y\times V\rightarrow V$, 
$(w,v)\mapsto w\circ_\xi v$
is defined by
\begin{equation*}
\begin{split}
w\circ_\xi v&=w(v),\qquad\,\,\, w\in W,\\
t_y\circ_\xi v&=v+\xi y,\qquad y\in Y.
\end{split}
\end{equation*}
Set $\widetilde{N}_\lambda=\{f\in C^\omega(V) \,\, | \,\, 
p(\partial)f=\chi_\lambda(p)f\,\,\,\forall p\in S(\hf)^W\}$ for $\lambda\in\hf$.
Here $\chi_\lambda$ is the linear character of $S(\hf)$ 
satisfying $v\mapsto (\lambda,v)$ and $p(\partial)$ stands for the
constant coefficient partial differential operator naturally associated
to $p\in S(\hf)$ by $v\mapsto \partial_v$ ($v\in V$). 
Then $\widetilde{N}_\lambda$ is a
$\#W$-dimensional vector space containing $e^{w\lambda}$ ($w\in W$).
In particular the $e^{w\lambda}$ ($w\in W$) form a
basis of $\widetilde{N}_\lambda$ 
if $\lambda\in\hf_{reg}:=\{\lambda\in\hf \, | \, (\lambda,\alpha)\not=0\,\,
\forall\alpha\in R\}$.

Since $S(\hf)^W\subseteq Z(\widetilde{H}_0^{Y}(k))$, 
\[
\widetilde{S}_M(\lambda):=\widetilde{N}_\lambda\otimes_{\mathbb{C}}M
\]
is a finite dimensional
$Q_\xi(\widetilde{H}_0^Y(k))$-submodule of $C^\omega(V)\otimes_{\mathbb{C}}M$
for all $\lambda\in\hf^*$.  

Define for $v\in V$ and
$\lambda\in\hf_{0,reg}$ the Bethe wave function
\[
\psi_\lambda(v):=\sum_{w\in W}e^{(w\lambda,v)}
J_w(\lambda)\in\mathbb{C}[W].
\]
Note that 
the $J_{t_y}(\lambda)$ ($y\in Y$) pairwise commute. 
The following result from \cite[Thm. 5.10]{EOS} 
relates the coordinate 
Bethe ansatz for vector valued root 
system analogs of the quantum Bose gas on the circle 
with pairwise delta-function interactions to the 
study of the space of $Q_\xi(\Wa^Y)$-invariants in $\widetilde{S}_M(\lambda)$.
\begin{thm}
Let $M$ be a finite dimensional $\Wa^Y$-module, $m\in M$ 
and $\lambda\in\hf_{0,reg}$.
Then $\psi_\lambda(\cdot)m\in 
\widetilde{S}_M(\lambda)^{Q_\xi(\Wa^Y)}$ if and only if
\[
J_{t_y}(\lambda)m=e^{\xi(\lambda,y)}m\qquad \forall y\in Y
\] 
(the Bethe ansatz equations).
\end{thm}

\subsection{Generalized Bethe wave functions}\label{Solsection}
Now we consider unspecified level. Without loss of generality we assume
that $Y$ spans $V$ (cf. Remark \ref{notfull}{\bf (iii)}), in which case
the center of $\wh{H}^Y(k)$ is generated by $c$ and $C$
(see Proposition \ref{centerprop}).
The role of $\widetilde{N}_\lambda$ is taken over by
\[
\widehat{N}_{\gamma,\Gamma}:=\{f\in C^\omega(\wh{V}^+) \,\, | \,\, 
\partial_c(f)=\gamma f\,\,\, \& \,\,\, \wh{\Delta}(f)=\Gamma f\}
\]
for $\gamma,\Gamma\in\mathbb{C}$,
where $\wh{\Delta}$ is the constant coefficient differential operator
associated to $C$,
\[
\wh{\Delta}=\Delta+2\partial_c\partial_d
\]
and $\Delta$ is the Laplacean on $V$. If $\wh{\lambda}\in\wh{\hf}_\kappa$
then 
\[
e^{w\wh{\lambda}}\in \wh{N}_{\kappa,(\wh{\lambda},\wh{\lambda})}\qquad 
\forall\, w\in\Wa^Y.
\]

Note that $\wh{S}_{M}(\gamma,\Gamma):=\wh{N}_{\gamma,\Gamma}
\otimes_{\mathbb{C}}M$ 
is a $\widehat{Q}(\wh{H}^Y(k))$-submodule of 
$C^\omega(\widehat{V}^+)\otimes M$. The generalized Bethe wave functions
will be $\wh{Q}(\Wa^Y)$-invariant elements in the infinite
dimensional vector space $\wh{S}_M(\gamma,\Gamma)$, explicitly defined
as a convergent series expansion in the
$e^{w\wh{\lambda}}$ ($w\in \Wa^Y$) for some 
$\wh{\lambda}\in\wh{\hf}_{\gamma}$ (see \eqref{hxi})
satisfying $(\wh{\lambda},\wh{\lambda})=\Gamma$.
In addition we show that the role of the Bethe ansatz equations is taken over
by the affine difference Knizhnik-Zamolodchikov equations.

If $M$ is a finite dimensional $\Wa^Y$-module 
then we regard $\textup{End}_{\mathbb{C}}(M)$ as $\Wa^Y$-module
by $w\cdot\psi:=\pi_M(w)\circ\psi$. Note that evaluation at $m\in M$ defines 
a $\wh{H}^Y(k)$-linear map $\wh{S}_{\textup{End}_{\mathbb{C}}(M)}(\gamma,\Gamma)
\rightarrow \wh{S}_M(\gamma,\Gamma)$.
\begin{thm}\label{E}
Let $M$ be a unitarizable finite dimensional $\Wa^Y$-module. 
For $\widehat{v}\in\widehat{V}^+$ and 
$\widehat{\lambda}\in\wh{\hf}_{reg}^+$ the series
\begin{equation*}
\begin{split}
E_M(\widehat{v};\widehat{\lambda})&:=\sum_{y\in Y}e^{(
t_y\widehat{\lambda},\wh{v})}J_{t_y}(\widehat{\lambda}),\\
E_M^+(\widehat{v};\widehat{\lambda})&:=
\sum_{w\in\Wa^Y}e^{(w\widehat{\lambda},\wh{v})}J_w(\widehat{\lambda})
\end{split}
\end{equation*}
converge in $\textup{End}_{\mathbb{C}}(M)$ and satisfy
\begin{enumerate}
\item $E_M^+(\widehat{v};\widehat{\lambda})=\sum_{w\in W}
E_M(\widehat{v};w\widehat{\lambda})J_w(\widehat{\lambda})$.
\item if $\widehat{\lambda}\in\wh{\hf}_{\kappa,reg}$ with $\kappa\in\mathbb{C}_+$
then
\begin{enumerate}
\item[(a)] $E_M(\cdot;\widehat{\lambda})\in C_\kappa^\omega(\widehat{V}^+)
\otimes_{\mathbb{C}}\textup{End}_{\mathbb{C}}(M)$,
\item[(b)] $E_M^+(\cdot;\widehat{\lambda})\in 
\wh{S}_{\textup{End}_{\mathbb{C}}(M)}(\kappa,
(\widehat{\lambda},\widehat{\lambda}))^{\wh{Q}(\Wa^Y)}$.
\end{enumerate}
\item if $\widehat{v}\in\widehat{V}_\xi$ with $\xi>0$ then
$E_M(\widehat{v};\cdot)\in\textup{adKZ}_\xi$
and $E_M^+(\widehat{v};\cdot)\in\textup{adKZ}_{\xi}^{\cdot W}$.
\end{enumerate}
\end{thm}
\begin{proof}
We first consider the convergence of the series, for which we use
\begin{lem}
Let $K\subset\wh{\hf}_{reg}^+$ be a compact subset. There exists
a positive constant $D=D(K)$ such that
\begin{equation}\label{estimate}
\|J_w(\widehat{\lambda})\|_M\leq D^{l(w)},\qquad \forall\, w\in\Wa^Y,\,
\forall \widehat{\lambda}\in\bigcup_{w\in\Wa^Y}wK,
\end{equation}
where $\|\cdot\|_M$ is the norm on $M$ turning $M$ into a unitary 
$\Wa^Y$-module.
\end{lem}
\begin{proof}
Since $\wh{\hf}_{reg}^+$ is $\Wa^Y$-invariant, $\cup_{w\in\Wa^Y}wK
\subset\wh{\hf}_{reg}^+$.
Note that there exists a $s=s(K)>0$ such that 
\[
|(\widehat{\lambda},a^\vee)-k_a|\geq s\quad \forall a\in\widehat{R},\,
\forall \widehat{\lambda}\in K.
\]
Then the estimate \eqref{estimate} is easily seen to be correct
for $w\in\Wa^Y$ with 
\[l(w):=\#(\wh{R}^+\cap w^{-1}\wh{R}^-)\leq 1
\]
(i.e. for $w\in\Wa^Y$ of the form $w=\omega$ or $w=s_a\omega$
with $a\in\wh{F}$ and $\omega\in\Omega$) if the constant $D=D(K)$ 
is taken to be
\[
D=1+\frac{2}{s}\max_{a\in\wh{R}}|k_a|.
\] 
Note that $D$ is well defined since 
$k$ attains only finitely many values.
The estimate \eqref{estimate} for arbitrary 
$w\in\Wa^Y$ holds true with the same positive constant $D$. This
follows 
by induction to $l(w)$ using the cocycle property of $J_w(\widehat{\lambda})$.
\end{proof}
Write $\lambda=\textup{Re}(\lambda)+\sqrt{-1}\textup{Im}(\lambda)$ with
$\textup{Re}(\lambda),\textup{Im}(\lambda)\in V$. 
Combining the lemma with the fact that
\begin{equation}\label{lengthty}
l(t_y)=\sum_{\alpha\in R^+}|(y,\alpha)|
\end{equation}
and
\begin{equation}\label{expexplicit}
e^{(t_y\wh{\lambda},\wh{v})}=
e^{u\kappa+\xi\eta+(\lambda,v)+(\kappa v-\xi\lambda,y)-\frac{\xi\kappa}{2}(y,y)}
\end{equation}
for $y\in Y$, $\wh{v}=v+uc+\xi d\in\wh{V}_\xi$ ($\xi>0$)
and $\wh{\lambda}=\lambda+\eta c+\kappa d\in\wh{\hf}_{reg}^+$, we
conclude that
\[
\sum_{y\in Y}\|e^{(t_y\wh{\lambda},\wh{v})}J_{t_y}(\wh{\lambda})\|_M
\leq e^{u\textup{Re}(\kappa)+\xi\textup{Re}(\eta)+\textup{Re}((\lambda,v))}
\sum_{y\in Y}e^{-\frac{\xi\textup{Re}(\kappa)}{2}(y,y)}
e^{(\textup{Re}(\kappa)v-\xi\textup{Re}(\lambda),y)+
\textup{log}(D)\sum_{\alpha\in R^+}|(y,\alpha)|}
\]
if $\wh{v}\in \wh{V}^+$ and 
$\wh{\lambda}\in K$ 
(with $K\subset\wh{\hf}_{reg}^+$ a fixed compact set and $D=D(K)$
the associated positive constant). The absolute convergence of
the series $E_M(\wh{v};\wh{\lambda})$ ($\wh{v}\in\wh{V}^+$ and 
$\wh{\lambda}\in\wh{\hf}_{reg}^+$) follows. It also implies that
$E_M(\cdot;\wh{\lambda})$ is real analytic on $\wh{V}^+$
for $\wh{\lambda}\in\wh{\hf}_{reg}^+$ and 
$E(\wh{v};\cdot)$ is holomorphic on $\wh{\hf}_{reg}^+$
for $\wh{v}\in\wh{V}^+$.

We continue now with the proof of (1)-(3).\\
(1) is immediate from the cocycle property of $\{J_w\}_{w\in\Wa^Y}$.\\
(2) We already established part (a).\\
(b) For $a\in\widehat{F}$,
\[\widehat{Q}(s_a)E_M^+(\cdot;\widehat{\lambda})=
\sum_{w\in\Wa^Y}e^{s_aw\widehat{\lambda}}s_aJ_w(\widehat{\lambda})
-k_a\sum_{w\in\Wa^Y}\bigl(I(a)e^{w\widehat{\lambda}}\bigr)
J_w(\widehat{\lambda}).
\]
Now use that for $\widehat{\lambda}\in\wh{\hf}_{reg}^+$,
\[
I(a)\bigl(e^{\widehat{\lambda}}\bigr)=\frac{e^{\widehat{\lambda}}-
e^{s_a\widehat{\lambda}}}{(a^\vee,\widehat{\lambda})},
\]
hence
\begin{equation}\label{intermediatestep}
\widehat{Q}(s_a)E_M^+(\cdot;\widehat{\lambda})=
\sum_{w\in\Wa^Y}e^{s_aw\widehat{\lambda}}s_aJ_w(\widehat{\lambda})-
k_a\sum_{w\in\Wa^Y}\frac{e^{w\widehat{\lambda}}-
e^{s_aw\widehat{\lambda}}}{(a^\vee,w\widehat{\lambda})}
J_w(\widehat{\lambda}).
\end{equation}
Since $\widehat{\lambda}\in\wh{\hf}_{reg}^+$, the explicit expression
of the cocycle value $J_{s_a}(w\widehat{\lambda})$ allows us to write
\[
s_a=\Bigl(1-\frac{k_a}{(a^\vee,
w\widehat{\lambda})}\Bigr)J_{s_a}(w\widehat{\lambda})
-\frac{k_a}{(a^\vee,w\widehat{\lambda})}.
\]
Substituting in \eqref{intermediatestep} and using the cocycle condition
we get
\[
\widehat{Q}(s_a)E_M^+(\cdot;\widehat{\lambda})=
\sum_{w\in\Wa^Y}\Bigl(1-\frac{k_a}{(a^\vee,w\widehat{\lambda})}\Bigr)
e^{s_aw\widehat{\lambda}}J_{s_aw}(\widehat{\lambda})-
\sum_{w\in\Wa^Y}\frac{k_a}{(a^\vee,w\widehat{\lambda})}e^{w\widehat{\lambda}}
J_w(\widehat{\lambda}).
\]
In the first sum replace the summation variable $w$ by $s_aw$.
It follows that
$\widehat{Q}(s_a)E_M^+(\cdot;\widehat{\lambda})=E^+(\cdot;\widehat{\lambda})$.
If $\omega\in\Omega$ then, using 
$\widehat{Q}(\omega)=\omega\otimes\pi_M(\omega)$ 
and $J_\omega(w\widehat{\lambda})=\omega$ for all $w\in\Wa^Y$,
\[
\widehat{Q}(\omega)E_M^+(\cdot;\widehat{\lambda})=
\sum_{w\in\Wa^Y}e^{\omega w\widehat{\lambda}}J_{\omega w}(\widehat{\lambda})=
E_M^+(\cdot;\widehat{\lambda}).
\]
Hence $E_M^+(\cdot;\widehat{\lambda})$ is $\widehat{Q}(\Wa^Y)$-invariant.\\
(3) The fact that $E_M(\widehat{v};\cdot)$ solves the affine difference 
KZ equations
is direct by the cocycle condition: if $y^\prime\in Y$ then
\[
E_M(\widehat{v};t_{y^\prime}\widehat{\lambda})J_{t_{y^\prime}}(\widehat{\lambda})=
\sum_{y\in Y}e^{(t_{y+y^\prime}\widehat{\lambda},\wh{v})}
J_{t_{y+y^\prime}}(\widehat{\lambda})=E(\widehat{v};\widehat{\lambda}).
\]
Similarly one establishes $E_M^+(\widehat{v};\cdot)\in\textup{adKZ}^{\cdot W}$.
\end{proof}
We call $E_M^+(\cdot,\wh{\lambda})\in 
\wh{S}_{\textup{End}_{\mathbb{C}}(M)}(\kappa,(\wh{\lambda},\wh{\lambda}))$
for $\wh{\lambda}\in\wh{\hf}_{\kappa,reg}$ ($\kappa\in\mathbb{C}_+$)
the {\it generalized Bethe wave functions of level} $\kappa$.
\begin{rema}\label{connection}
Formally one can think of the affine difference KZ equations as defining
a difference connection on the bundle over $\wh{\hf}_{reg}^+$ with
fiber at $\wh{\lambda}\in\wh{\hf}_{\kappa,reg}$ ($\kappa\in\mathbb{C}_+$)
given by $\wh{S}_M(\kappa,(\wh{\lambda},\wh{\lambda}))$. The difference
connection commutes with
the action of $\wh{Q}(\wh{H}^Y(k))$ on the fibers of the bundle. 
The generalized Bethe
wave function $\wh{\lambda}\mapsto E_M^+(\cdot,\wh{\lambda})$ then
defines a flat $\wh{Q}(\Wa^Y)$-invariant section.
\end{rema}

\begin{rema}\label{finitelimit}
Fix $\kappa\in\mathbb{C}_+$. If $\lambda\in\hf_{\kappa,reg}$ satisfies
\begin{equation}\label{convergencedomain}
(\textup{Re}(\lambda),y)>-\frac{\textup{Re}(\kappa)}{2}(y,y)\quad \forall\, 
y\in Y\setminus\{0\},
\end{equation}
then 
\begin{equation*}
\begin{split}
\lim_{\xi\rightarrow\infty}E_M(v+\xi d,\lambda+\kappa d)&
=e^{(\lambda,v)}\textup{Id}_M,\\
\lim_{\xi\rightarrow\infty}E_M^+(v+\xi d,\lambda+\kappa d)&=
\sum_{w\in W}e^{(w\lambda,v)}J_w(\lambda)=\psi_\lambda(v)
\end{split}
\end{equation*}
uniformly for $v$ in compacta of $V$. Note that
$v\mapsto\psi_\lambda(v)$
coincides with the Bethe wave function of the vector valued
root system analog of the quantum Bose gas on the line
with pairwise delta-function interactions, cf., e.g., \cite{HO},
\cite{EOS} and references therein. The fundamental property of
$\psi_\lambda\in 
C^\omega(V)\otimes_{\mathbb{C}}\textup{End}_{\mathbb{C}}(M)$
is the fact that it is a $Q(W)$-invariant real analytic solution
to the differential equations $p(\partial)f=p(\lambda)f$ for all 
$p\in S(\hf)^W$, where $Q$ is the $W$-action on
$C^\omega(V)\otimes_{\mathbb{C}}\textup{End}_{\mathbb{C}}(M)$ given by
integral-reflection operators, cf. \cite{HO,EOS0,EOS}.
\end{rema}

\begin{eg}\label{Stein}
The simplest example corresponds to the Steinberg module
$M=\textup{St}$, which is the
one-dimensional $\Wa^Y$-module with associated
linear character $w\mapsto (-1)^{l(w)}$ ($w\in\Wa^Y$), since
$J_w(\widehat{\lambda})|_{\textup{St}}=(-1)^{l(w)}$ ($w\in\Wa^Y$). 
We regard $E_{\textup{St}}$ and $E_{\textup{St}}^+$ 
as scalar valued functions. 
Assume for simplicity that $(Y,P)\subseteq\mathbb{Z}$, where $P$ is the
weight lattice of $R$. Then $l(t_y)$ is even for all 
$y\in Y$ (indeed, by $W$-invariance it suffices to prove it when
$(y,\alpha)\geq 0$ for all $\alpha\in R$, in which case 
$l(t_y)=2(\rho,y)$ with $\rho\in P$ the half sum of positive
roots by \eqref{lengthty}). Then
\begin{equation*}
\begin{split}
E_{\textup{St}}(\widehat{v},\widehat{\lambda})&=\sum_{y\in Y}
e^{(t_y\widehat{\lambda},\widehat{v})}\\
&=e^{u\kappa+\xi\eta+(\lambda,v)}\sum_{y\in Y}e^{\kappa (v,y)-\xi (\lambda,y)
-\frac{\xi\kappa}{2}(y,y)}  
\end{split}
\end{equation*}
if $\widehat{v}=v+uc+\xi d\in\wh{V}_\xi$ ($\xi>0$) and
$\wh{\lambda}=\lambda+\eta c+\kappa d\in\wh{\hf}_{\kappa,reg}$ 
($\kappa\in\mathbb{C}_+$). This is essentially 
a classical theta function (cf. \cite[Chpt. 13]{Kac1990}
and references therein).
Furthermore,
\[
E_{\textup{St}}^+(\widehat{v},\widehat{\lambda})=
\sum_{w\in W}(-1)^{l(w)}E_{\textup{St}}(\wh{v},w\wh{\lambda}).
\]
\end{eg}
\begin{eg}
Consider the trivial $\Wa^Y$-module 
$M=\textup{Triv}$, which is the one-dimensional $\Wa^Y$-module
with associated linear character $w\mapsto 1$ ($w\in\Wa^Y$).
Let $j_w\in\mathcal{O}(\wh{\hf}_{reg}^+)$ ($w\in\Wa^Y$) such that
$J_w(\wh{\lambda})$ acts as multiplication by $j_w(\wh{\lambda})$ on
$\textup{Triv}$. Viewing $E_{\textup{Triv}}$ and $E_{\textup{Triv}}^+$ 
as scalar valued functions we thus get
\begin{equation*}
\begin{split}
E_{\textup{Triv}}(\wh{v},\wh{\lambda})&=
\sum_{y\in Y}j_{t_y}(\wh{\lambda})e^{(t_y\wh{\lambda},\wh{v})},\\
E_{\textup{Triv}}^+(\wh{v},\wh{\lambda})&=\sum_{w\in W}E_{\textup{Triv}}(\wh{v},
w\wh{\lambda})j_w(\wh{\lambda}).
\end{split}
\end{equation*}
The cocycle values are explicitly given by
\[
j_w(\wh{\lambda})=\prod_{a\in\wh{R}^+\cap w^{-1}\wh{R}^-}
\frac{(a^\vee,\wh{\lambda})+k_a}{(a^\vee,\wh{\lambda})-k_a},\qquad
w\in\Wa^Y.
\]
In particular, if 
$k_{\alpha+mc}=k_\alpha$ for all $\alpha\in R$ and $m\in\mathbb{Z}$
(this is automatically true 
if $P^\vee\subseteq Y$, with $P^\vee$ the coweight lattice
of $R$) then
\[
j_{t_y}(\wh{\lambda})=\prod_{\alpha\in R^+: (\alpha,y)>0}\prod_{m=0}^{(\alpha,y)-1}
\frac{m\kappa_\alpha^\vee+
(\alpha^\vee,\lambda)+k_\alpha}{m\kappa_\alpha^\vee+(\alpha^\vee,\lambda)
-k_\alpha}
\prod_{\beta\in R^+: (\beta,y)<0}\prod_{m=1}^{-(\beta,y)}\frac{m\kappa_\beta^\vee-
(\beta^\vee,\lambda)+k_\beta}
{m\kappa_\beta^\vee-(\beta^\vee,\lambda)-k_\beta},
\]
where $\wh{\lambda}=\lambda+\eta c+\kappa d$ and
$\kappa_\alpha^\vee:=2\kappa/(\alpha,\alpha)$.
\end{eg}

We present now a straightforward
generalization of some of the statements of Theorem \ref{E}.
Call $f\in\mathcal{O}(\wh{\hf}_{reg}^+)\otimes_{\mathbb{C}}
\textup{End}_{\mathbb{C}}(M)$ $\Wa^Y$-invariant if $f\lhd w=f$ for
$w\in \Wa^Y$, where $\lhd$ is the right $\Wa^Y$-action on
$\mathcal{O}(\wh{\hf}_{reg}^+)\otimes_{\mathbb{C}}
\textup{End}_{\mathbb{C}}(M)$ defined by
\[
(f\lhd w)(\wh{\lambda}):=J_w(\wh{\lambda})^{-1}f(w\wh{\lambda})
J_w(\wh{\lambda}).
\]
\begin{cor}
Let $M$ be a finite dimensional unitarizable $\Wa^Y$-module
and take $f\in\mathcal{O}(\wh{\hf}_{reg}^+)\otimes_{\mathbb{C}}
\textup{End}_{\mathbb{C}}(M)$. Set
\[
E_{M,f}^+(\wh{v},\wh{\lambda}):=E_M^+(\wh{v},\wh{\lambda})f(\wh{\lambda})
\]
for $\wh{v}\in\wh{V}^+$ and $\wh{\lambda}\in\wh{\hf}_{reg}^+$.

For $\wh{\lambda}\in\wh{\hf}_{\kappa,reg}$ and $\wh{v}\in\wh{V}_\xi^+$
($\kappa\in\mathbb{C}_+$ and $\xi\in\mathbb{R}_{>0}$) we then have
\begin{enumerate}
\item $E_{M,f}^+(\cdot,\wh{\lambda})\in 
\wh{S}_{\textup{End}_{\mathbb{C}}(M)}(\kappa,(\wh{\lambda},
\wh{\lambda}))^{\wh{Q}(\Wa^Y)}$,
\item if $f$ is $\Wa^Y$-invariant then
$E_{M,f}^+(\wh{v},\cdot)\in\textup{adKZ}_\xi^{\cdot W}$.
\end{enumerate}
\end{cor}
\begin{proof}
(1) is clear from the theorem.\\
(2) Since $f$ is $\Wa^Y$-invariant,
\[
E_{M,f}^+(\wh{v},\wh{\lambda})=\sum_{w\in\Wa^Y}e^{(w\wh{\lambda},\wh{v})}
f(w\wh{\lambda})J_w(\wh{\lambda}),
\]
which is $\cdot\Wa^Y$-invariant.
\end{proof}
Let $\wh{\lambda}: (0,1)\rightarrow\wh{\hf}_{reg}^+$ be a path
such that $\wh{\lambda}(t)\rightarrow \lambda\in\hf_{0,reg}$
if $t\downarrow 0$. Then for all $w\in \Wa^Y$,
\[
\lim_{t\downarrow 0}e^{(\wh{\lambda}(t),w(v+\xi d))}=
e^{(\lambda,w\circ_\xi v)}
\]
and for all $a\in \Ra$,
\[
\lim_{t\downarrow 0}(I(a)e^{\wh{\lambda}(t)})(v+\xi d)=
\int_0^{a^\xi(v)}e^{(\lambda,v-tDa^\vee)}dt.
\]
Consequently
\[
\lim_{t\downarrow 0}\bigl(\wh{Q}(w)e^{\wh{\lambda}(t)}\bigr)(v+\xi d)=
\bigl(Q_\xi(w)e^{\lambda}\bigr)(v),\qquad w\in \Wa^Y.
\]
This suggests that 
if $E_M^+(v+\xi d,\wh{\lambda}(t))$ converges as $t\downarrow 0$,
then it converges to the Bethe wave function $\psi_\lambda(v)$
associated to the vector valued root system
analog of the quantum Bose gas on the circle with pairwise
delta-function interactions.
More concretely we expect that for $\wh{v}\in\wh{V}_\xi$ and for 
$(\lambda,m)\in\hf_{0,reg}\times M$ satisfying the Bethe ansatz equations
$J_{t_y}(\lambda)m=e^{\xi(\lambda,y)}m$ for all $y\in Y$,
a renormalization of $E_M(\wh{v},\wh{\lambda})m$
converges to $e^{\xi(\lambda,v)}m$ when $\wh{\lambda}$ tends to 
$\lambda$ along specific paths.

Here we will only consider the simplest case
that $M=\textup{St}$.
For the sake of simplicity we assume that $(Y,P)\subseteq\mathbb{Z}$
(cf. Example \ref{Stein}).
Note that the $k$-dependence drops out for $M=\textup{St}$
(it relates to 
the limit $k\rightarrow\infty$ of the theory for the trivial representation
$M=\textup{Triv}$, cf. \cite[\S 3]{EOS0}).
The $J_{t_y}(\wh{\lambda})$ ($y\in Y$)
act trivially on $\textup{St}$. Hence the Bethe ansatz equations
simplify to the requirement that $\lambda\in\hf_{reg}$ satisfies
\[
e^{\xi(\lambda,y)}=1\qquad \forall y\in Y.
\]
Equivalently $\lambda\in \frac{2\pi\sqrt{-1}}{\xi}X_{reg}$, where
$X_{reg}=X\cap\hf_{reg}$ and
$X\subset V$ is the lattice dual to $Y$ with respect to $(\cdot,\cdot)$.
Fix $\lambda\in\frac{2\pi\sqrt{-1}}{\xi}X_{reg}$ and $\wh{v}\in\wh{V}_\xi$,
say $\wh{v}=v+\eta c+\xi d$. Then 
\[
E_{St}(\wh{v},\lambda+\kappa d)=e^{(\lambda,v)}E_{St}(\wh{v},\kappa d),
\]
hence for all $y\in Y$,
\begin{equation*}
\begin{split}
E_{St}(\wh{v},\kappa d)^{-1}E_{St}(\wh{v},\lambda+\kappa y+\kappa d)
&=E_{St}(\wh{v},\kappa)^{-1}E_{St}(\wh{v},t_y(\lambda+\kappa d))
e^{\frac{\kappa\xi}{2}(y,y)+\xi(y,\lambda)}\\
&=E_{St}(\wh{v},\kappa)^{-1}E_{St}(\wh{v},\lambda+\kappa d)
e^{\frac{\kappa\xi}{2}(y,y)+\xi(y,\lambda)}\\
&=e^{\frac{\kappa\xi}{2}(y,y)+(\lambda,v)}
\end{split}
\end{equation*}
which converges to $e^{(\lambda,v)}$ as $\kappa\rightarrow 0$.
In particular, for $\wh{v}=v+\eta c+\xi d\in \wh{V}_\xi$ and
$\lambda\in\frac{2\pi\sqrt{-1}}{\xi}X_{reg}$,
\[
\lim_{\mathbb{C}_+\ni\kappa\rightarrow 0}
E_{St}(\wh{v},\kappa d)^{-1}E_{St}^+(\wh{v},\lambda+\kappa d)=
\sum_{w\in W}(-1)^{l(w)}e^{(w\lambda,v)},
\]
in accordance to \cite[Prop. 3.3]{EOS0}.

\section{Time dependent Schr{\"o}dinger equations 
with delta potentials}
In this section we construct the nonstationary Schr{\"o}dinger equation
with delta potentials and its solutions. We first need to describe
the image of the propagation operator $T$ in more detail.
For $b\in\wh{R}$ set $H_b:=\{\wh{v}\in\wh{V}^+\, | \, (\wh{v},b)=0\}$.
\begin{prop} \label{prp:jumpconditions}
Let $M$ be a $\Wa^Y$-module and
let $w\in \Wa^Y$,
 $b\in w \Fa$, $\widehat{v}\in H_b$, $p\in S(\widehat{\hf})$ and
$f\in T\big(C^\om(\widehat{V})\otimes_\C M\big)$.
Then
\begin{equation}\label{jump} 
\big( p(\partial) f_{w\widehat{C}^+}\big)(\widehat{v}) -
\big( p(\partial) f_{s_bw\widehat{C}^+}\big)(\widehat{v}) =
k_b\pi_M(s_b)\Big(\big(\De_b(p)(\partial) f_{w\widehat{C}^+}\big) 
(\widehat{v})\Big),
\end{equation}
where $\Delta_b(p):=\frac{s_bp-p}{b^\vee}\in S(\wh{\hf})$. In particular
$f_{w\wh{C}^+}(\wh{v})=f_{s_bw\wh{C}^+}(\wh{v})$ and
if $\wh{u}\in\wh{V}$, regarded as an element in $S(\wh{\hf})$ of degree one,
then
\[
(\partial_{\wh{u}}f_{w\wh{C}^+})(\wh{v})-(\partial_{\wh{u}}f_{s_bw\wh{C}^+})(\wh{v})
=-k_b(\wh{u},b)\pi_M(s_b)f_{w\wh{C}^+}(\wh{v}).
\]
\end{prop}
\begin{proof} (compare with the 
proof of \cite[Thm 5.3(ii)]{EOS0} and \cite[Prop. 4.16]{EOS}).

By $\Wa^Y$-equivariance of $T$ and \eqref{actionF}, it suffices
to prove \eqref{jump} for $w=1$.

Let $a\in \Fa$, $\widehat{v}\in H_a$,
$p\in S(\wh{\hf})$ and $f=Tg$ with
$g\in C^\om(\widehat{V})\otimes_\C M$. In $\wh{H}^Y(k)$ we have the cross
relation
\[
s_a\cdot p-s_a(p)\cdot s_a=k_a\Delta_a(p).
\]
Applying the representation map $\wh{Q}$, acting on $g$ and evaluating at
$\wh{v}$ we get 
\[
\pi_M(s_a)(p(\partial)g)(\wh{v})-\bigl((s_ap)(\partial)\wh{Q}(s_a)g
\bigr)(\wh{v})=k_a\bigl(\Delta_a(p)(\partial)g\bigr)(\wh{v})
\]
because $(\wh{Q}(s_a)p(\partial)g)(\wh{v})=\pi_M(s_a)(p(\partial)g)(\wh{v})$.
Acting on both sides by $\pi_M(s_a)$ we get
\[
(p(\partial)g)(\wh{v})-
\bigl(p(\partial)\pi_M(s_a)(\wh{Q}(s_a)g)(s_a\cdot)\bigr)(\wh{v})=
k_a\pi_M(s_a)\bigl(\Delta_a(p)(\partial)g\bigr)(\wh{v}).
\]
Since $f_{\wh{C}^+}=g$ and $f_{s_a\wh{C}^+}=\pi_M(s_a)(\wh{Q}(s_a)g)(s_a\cdot)$,
we get the desired identity
\[
\bigl(p(\partial)f_{\wh{C}^+}\bigr)(\wh{v})-
\bigl(p(\partial)f_{s_a\wh{C}^+}\bigr)(\wh{v})=
k_a\pi_M(s_a)\bigl(\Delta_a(p)(\partial)f_{\wh{C}^+}\bigr)(\wh{v}).
\]
\end{proof}
For $f\in\mathcal{F}_M$ write $f^\dagger\in 
L^1_{loc}(\wh{V}^+)\otimes_{\mathbb{C}}M$ such that 
$f^\dagger(\wh{v})=f_{\wh{C}}(\wh{v})$ if $\wh{v}\in\wh{C}$ 
($\wh{C}\in\mathcal{C}$). Write
\[
C_M:=\{f\in\mathcal{F}_M \,\, | \,\, f_{w\wh{C}^+}(\wh{v})=
f_{s_aw\wh{C}^+}(\wh{v})\quad \forall\, w\in\Wa^Y,\, a\in w\wh{F}\, 
\hbox{ and } \wh{v}\in H_a\}.
\]
If $f\in C_M$ then the local $L^1$-function $f^\dagger$ can be represented
by a continuous function on $\wh{V}^+$, which we will also denote
by $f^\dagger$. It is explicitly given
by $f^\dagger(\wh{v})=f_{\wh{C}}(\wh{v})$
if $\wh{v}\in\wh{V}^+\cap \overline{\wh{C}}$. 
Note that $T(C^\omega(\wh{V}^+)\otimes_{\mathbb{C}}M)
\subseteq C_M$ in view of Proposition \ref{prp:jumpconditions}. 

We consider $\wh{V}$ as Euclidean space with scalar product
\[
\langle \wh{u},\wh{v}\rangle:=\bigl(E(\wh{u}),\wh{v}\bigr),
\]
where $E$ is the linear involution of $\wh{V}$ satisfying
$E|_V=\textup{Id}_V$, $E(c)=d$ and $E(d)=c$. If $\{v_1,\ldots,v_m\}$
is an orthonormal basis of $V$ with respect to $(\cdot,\cdot)$ then
$\{v_1,\ldots,v_m,c,d\}$ is an orthonormal basis of $\wh{V}$ with respect
to $\langle \cdot,\cdot\rangle$.
Let $\mathrm{d}\widehat{v}$ denote the Lebesgue measure on 
$\widehat{V}$ and $\mathrm{d}\sigma(\wh{v})$ the 
induced measure on a hypersurface.

Define
\[\iota: L^1_\mathrm{loc}(\widehat{V}^+)\otimes_{\mathbb{C}}M \to
 \Hom_\C\big(C^\infty_\mathrm{c}(\widehat{V}^+),M\big) \]
by
\[(\iota g)(\ph) = \int_{\widehat{V}^+}g(\widehat{v})
\ph(\widehat{v}) \mathrm{\,d}\widehat{v}.
\]
Write
\[\widehat{\mathcal{H}}_k^M  = \widehat{\De} + \sum_{a\in\wh{R}^+}
 k_a\sqrt{\langle a,a\rangle}\delta\bigl((a,\cdot)\bigr)\pi_M(s_a)\]
for the linear map $\widehat{\mathcal{H}}_k^M :
C(\widehat{V}^+)\otimes_{\mathbb{C}}M \to
 \Hom_\C\big(C^\infty_\mathrm{c}(\widehat{V}^+),M\big)$ defined by
\[
(\widehat{\mathcal{H}}_k^M g )(\varphi) = 
\int_{\widehat{V}^+}  g(\wh{v})(\widehat{\Delta}\varphi)(\wh{v})
\mathrm{\,d}\widehat{v}
+\sum_{a\in \wh{R}^{+}} k_a\sqrt{\langle a,a\rangle}
  \int_{H_a} \pi_M(s_a)g(\widehat{v})\varphi(\widehat{v}) 
\mathrm{\,d}\sigma(\wh{v}).
\]
\begin{prop}\label{Hamiltonian}
Let $M$ be a finite dimensional $\Wa^Y$-module.
Assume $f\in C_M$ satisfies the 
derivative jump conditions
\begin{equation}\label{eq:jumpcondition1}
(\partial_{\widehat{u}}f_C)(\widehat{v})-
(\partial_{\widehat{u}}f_{s_bC})(\widehat{v})=
-k_b(b,\widehat{u})\pi_M(s_b)f_C(\widehat{v})
\end{equation}
for all $b\in w \Fa$, $C=w\widehat{C}_+$, $w\in \Wa$,
$\widehat{u}\in\widehat{V}$ and $\wh{v}\in H_b$.
Then 
\[\iota\bigl((\widehat{\De}f)^\dagger\bigr)=
\widehat{\mathcal{H}}_k^M(f^\dagger).\]
\end{prop}
\begin{rema}
By Proposition \ref{prp:jumpconditions},
the jump conditions \eqref{eq:jumpcondition1} are
satisfied if $f$ lies in the image of the propagation operator $T$.
\end{rema}
\begin{proof}
Let $\varphi\in C_c^\infty(\wh{V}^+)$ be a test function, then
\[
\iota\bigl((\widehat{\De}f)^\dagger\bigr)(\varphi)=
\sum_{\wh{C}\in\mathcal{C}}\int_{\wh{C}}(\wh{\Delta}f_{\wh{C}})(\wh{v})
\varphi(\wh{v})d\wh{v}.
\]
Substituting the definition of $\wh{\Delta}$ and applying 
repeatedly the divergence theorem one obtains
\begin{equation}\label{almost}
\begin{split}
\iota\bigl((\widehat{\De}f)^\dagger\bigr)(\varphi)&=
\int_{\wh{V}^+}f^{\dagger}(\wh{v})(\wh{\Delta}\varphi)(\wh{v})d\wh{v}\\
+\sum_{\wh{C}\in\mathcal{C}}\int_{\partial\wh{C}}
&\Bigl(\langle(\textup{grad}f_{\wh{C}})(\wh{v}),E(N^{\wh{C}}(\wh{v}))\rangle
\varphi(\wh{v})-f_{\wh{C}}(\wh{v})\langle(\textup{grad}\varphi)(\wh{v}),
E(N^{\wh{C}}(\wh{v}))\rangle\Bigr)d\sigma(\wh{v}),
\end{split}
\end{equation}
where $N^{\wh{C}}:\partial\wh{C}\rightarrow\wh{V}$ is the unit 
outward normal vector field on the boundary $\partial\wh{C}$ of $\wh{C}$.

We simplify now the boundary terms. Write
\begin{equation*}
\begin{split}
&\sum_{\wh{C}\in\mathcal{C}}\int_{\partial\wh{C}}
f_{\wh{C}}(\wh{v})\langle(\textup{grad}\varphi)(\wh{v}),
E(N^{\wh{C}}(\wh{v}))\rangle d\sigma(\wh{v})=\\
&\,\,\,=\sum_{a\in\Ra^+}\sum_{\stackrel{\wh{C}\in\mathcal{C}:}
{(a,\wh{C})\subseteq\mathbb{R}_{>0}}}
\int_{\partial\wh{C}^+\cap H_a}
\Bigl(f_{\wh{C}}(\wh{v})\langle(\textup{grad}\varphi)(\wh{v}),
E(N^{\wh{C}}(\wh{v}))\rangle\Bigr.\\
&\Bigl.\qquad\qquad\qquad\qquad\qquad\qquad\quad +
f_{s_a\wh{C}}(\wh{v})\langle(\textup{grad}\varphi)(\wh{v}),
E(N^{s_a\wh{C}}(\wh{v}))\rangle\Bigr)d\sigma(\wh{v}).
\end{split}
\end{equation*}
The contribution of the integral over $\wh{C}=w\wh{C}^+$ ($w\in\Wa$) 
will be zero unless $a\in w\wh{F}$, in which case
$N^{\wh{C}}(\wh{v})=-\frac{E(a)}{\sqrt{\langle a,a\rangle}}$ 
and $N^{s_a\wh{C}}(\wh{v})=\frac{E(a)}{\sqrt{\langle a,a\rangle}}$
for $\wh{v}\in \partial\wh{C}^+\cap H_a$. Since $f\in C_M$
we have in addition $f_{\wh{C}}(\wh{v})=f_{s_a\wh{C}}(\wh{V})$
for $\wh{v}\in H_a$, hence 
\[
\sum_{\wh{C}\in\mathcal{C}}\int_{\partial\wh{C}}
f_{\wh{C}}(\wh{v})\langle(\textup{grad}\varphi)(\wh{v}),
E(N^{\wh{C}}(\wh{v}))\rangle d\sigma(\wh{v})=0.
\]
Similarly
\begin{equation*}
\begin{split}
\sum_{\wh{C}\in\mathcal{C}}&\int_{\partial\wh{C}}
\langle(\textup{grad}f_{\wh{C}})(\wh{v}),E(N^{\wh{C}}(\wh{v}))\rangle
\varphi(\wh{v})d\sigma(\wh{v})\\
&=-\sum_{a\in\Ra^+}\frac{1}{\sqrt{\langle a,a\rangle}}
\sum_{\stackrel{\wh{C}\in\mathcal{C}:}
{(a,\wh{C})\subseteq\mathbb{R}_{>0}}}
\int_{\partial\wh{C}^+\cap H_a}
\Bigl(\frac{\partial f_{\wh{C}}}{\partial E(a)}(\wh{v})-
\frac{\partial f_{s_a\wh{C}}}{\partial E(a)}(\wh{v})\Bigr)
\varphi(\wh{v})d\sigma(\wh{v}).
\end{split}
\end{equation*}
By \eqref{eq:jumpcondition1} it follows that
\[
\sum_{\wh{C}\in\mathcal{C}}\int_{\partial\wh{C}}
\langle(\textup{grad}f_{\wh{C}})(\wh{v}),E(N^{\wh{C}}(\wh{v}))\rangle
\varphi(\wh{v})d\sigma(\wh{v})=
\sum_{a\in\wh{R}^+}k_b\sqrt{\langle a,a\rangle}\pi_M(s_a)
\int_{H_a}f^\dagger(\wh{v})\varphi(\wh{v})d\sigma(\wh{v}).
\]
Substitution in \eqref{almost} shows that
$\iota((\wh{\Delta}f)^\dagger)=\wh{\mathcal{H}}_k^M(f^\dagger)$.
\end{proof}

For $\wh{\lambda}\in\wh{\hf}_{reg}^+$ and a unitarizable finite dimensional
$\Wa^Y$-module $M$ set
\[
\phi_{\wh{\lambda}}:=T\bigl(E_M^+(\cdot;\wh{\lambda})\bigr)\in
C_{\textup{End}_{\mathbb{C}}(M)}.
\]
Note that $\phi_{\wh{\lambda}}$ is the unique $\Wa^Y$-invariant 
element in $\mathcal{F}_{\textup{End}_{\mathbb{C}}(M)}$ satisfying
\[
\phi_{\wh{\lambda},\wh{C}^+}=E_M^+(\cdot;\wh{\lambda}).
\]
Note that $\wh{\lambda}\mapsto\phi_{\wh{\lambda}}(\wh{v})$ is 
in $\textup{adKZ}_\xi^W$ if $\wh{v}\in\wh{V}_\xi$ ($\xi>0$).
\begin{prop}
With the above assumptions, 
\[
\widehat{\mathcal{H}}_k^{\textup{End}_{\mathbb{C}}(M)}(\phi_{\wh{\lambda}}^\dagger)=
(\wh{\lambda},\wh{\lambda})\iota(\phi_{\wh{\lambda}}^\dagger).
\]
In particular, if $\wh{\lambda}\in\wh{\hf}_{reg,\kappa}$ ($\kappa\in\mathbb{C}_+$)
then 
\begin{equation}\label{timedepHam}
\bigl(2\kappa\partial_d+\Delta+
\sum_{a\in\wh{R}^+}k_a\sqrt{\langle a,a\rangle}\delta((a,\cdot))\pi_M(s_a)\bigr)
\phi_{\wh{\lambda}}^\dagger=
(\wh{\lambda},\wh{\lambda})\iota(\phi_{\wh{\lambda}}^\dagger)
\end{equation}
in the weak sense, meaning that the left hand side tested against
$\varphi\in C_c^\infty(\wh{V}^+)$ is
\[
\int_{\wh{V}^+}\phi_{\wh{\lambda}}^\dagger(\wh{v})
((-2\kappa\partial_d+\Delta)\varphi)(\wh{v})d\wh{v}+
\sum_{a\in\wh{R}^+}k_a\sqrt{\langle a,a\rangle}\int_{H_a}
\pi_M(s_a)\phi_{\wh{\lambda}}^\dagger(\wh{v})\varphi(\wh{v})d\sigma(\wh{v}).
\]
\end{prop}
\begin{proof}
$\phi_{\wh{\lambda}}\in C_{\textup{End}_{\mathbb{C}}(M)}$ satisfies
the jump conditions \eqref{eq:jumpcondition1} by
Proposition \ref{prp:jumpconditions}. Hence 
\[
\widehat{\mathcal{H}}_k^M(\phi_{\wh{\lambda}}^\dagger)=
\iota((\wh{\Delta}\phi_{\wh{\lambda}})^\dagger)=
(\wh{\lambda},\wh{\lambda})\iota(\phi_{\wh{\lambda}})
\]
by Proposition \ref{Hamiltonian}. The second statement follows
from the first since $\partial_c(E_M^+(\cdot;\wh{\lambda}))=
\kappa E_M^+(\cdot;\wh{\lambda})$.
\end{proof}

Formula \eqref{timedepHam} can be interpreted as
$\phi_{\wh{\lambda}}$, with $\wh{\lambda}=\lambda+\eta c+\kappa d\in
\wh{\hf}_{reg,\kappa}$ and $\kappa\in\mathbb{C}_+$, 
solving (weakly) the time dependent Schr{\"o}dinger equation
\[
\partial_df=H_\kappa f
\]
with quantum Hamiltonian
\[
H_\kappa=\frac{1}{2\kappa}\Bigl(-\Delta-
\sum_{a\in\wh{R}^+}k_a\sqrt{\langle a,a\rangle}
\delta((a,\cdot))\pi_M(s_a)+(\lambda,\lambda)\Bigr)
+\eta.
\]
\begin{rema}
At $\kappa=0$, \eqref{timedepHam} formally reduces to 
a weak eigenvalue equation for 
\[
\Delta+\sum_{a\in\wh{R}^+}k_a\sqrt{\langle a,a\rangle}\delta((a,\cdot))\pi_M(s_a).
\]
This is {\it not} the weak quantum Hamiltonian of the vector valued
root system analog of the quantum Bose gas on the circle with pairwise
delta-function interactions
as analyzed in \cite[\S 4.5]{EOS}
with the help of the trigonometric Cherednik algebra $\widetilde{H}_0^{Y}(k)$,
which is given by
\[
\Delta+\sum_{a\in\wh{R}^+}k_a\sqrt{(a,a)}\delta((a,\cdot))\pi_M(s_a),
\]
see \cite[(4.18)]{EOS}.
\end{rema}
\begin{rema}
Take $\wh{\lambda}=\lambda+\kappa d\in\wh{\hf}_{reg,\kappa}$ 
with $\lambda$ satisfying 
\eqref{convergencedomain}. Then
\[
\psi_{\lambda}^0(v):=\lim_{\xi\rightarrow \infty}
\phi_{\lambda+\kappa d}^\dagger(v+\xi d),\qquad
v\in V
\]
is well defined and independent of $\kappa\in\mathbb{C}_+$. 
An explicit expression is given as follows.
Let
\[
K^+:=\{v\in V \,\, | \,\, (v,\alpha)>0 \quad \forall\, \alpha\in F\}
\]
be the fundamental Weyl chamber of $V$. Choose
$w\in W$ such that $v\in\overline{wK^+}$. Then $v+\xi d\in\overline{w\wh{C}^+}$
if $\xi\gg 0$ and
\[
\psi_{\lambda}^0(v)=\pi_M(w)\bigl(Q(w^{-1})\psi_\lambda\bigr)(w^{-1}v),
\]
cf. Remark \ref{finitelimit}.
This shows that $\psi_\lambda^0$ is 
the Bethe wave function for the vector valued root system analog
of the quantum Bose gas on the line with pairwise
delta-function interactions. 
It is continuous
and satisfies
\[
\bigl(\Delta+\sum_{\alpha\in R^+}k_\alpha\sqrt{(\alpha,\alpha)}
\delta((\alpha,\cdot))\pi_M(s_\alpha)\bigr)\psi_\lambda=
(\lambda,\lambda)\psi_\lambda
\]
in the weak sense. These observations should be compared to the analysis
at the trigonometric level in \cite[\S 5]{EK1}. 
\end{rema}

\end{document}